\documentclass[]{amsart}

\usepackage{amsbsy,amsfonts,amsmath,amssymb,enumerate,epsfig,graphicx,rotating,subfigure}
\usepackage[numbers,sort&compress]{natbib}

\usepackage{epsfig}
\usepackage{bm}

\usepackage{color}
\usepackage{algorithmicx}
\usepackage{algpseudocode}
\usepackage{mathtools}
\usepackage{url}
\usepackage[T1]{fontenc}

\newcommand{\sfrac}[2]{\mathchoice
  {\kern0em\raise.5ex\hbox{\the\scriptfont0 #1}\kern-.15em/
   \kern-.15em\lower.25ex\hbox{\the\scriptfont0 #2}}
  {\kern0em\raise.5ex\hbox{\the\scriptfont0 #1}\kern-.15em/
   \kern-.15em\lower.25ex\hbox{\the\scriptfont0 #2}}
  {\kern0em\raise.5ex\hbox{\the\scriptscriptfont0 #1}\kern-.2em/
   \kern-.15em\lower.25ex\hbox{\the\scriptscriptfont0 #2}}
  {#1\!/#2}}

\def\half   {\frac{1}{2}}
\def\myhalf {\sfrac{1}{2}}

\def\Yb {{\bf Y}}

\begin{document}

\markboth{Pazner et al.}{High-Order SDC for Low Mach Number Combustion}
\title[High-Order SDC for Low Mach Number Combustion]{A high-order spectral deferred correction strategy for low Mach number flow with complex chemistry}

\author[Pazner et al.]{Will Pazner$^{\rm \MakeLowercase{a} \ast}$}
\thanks{$^{\rm \ast}$Corresponding author. Email addresses: will\_pazner@brown.edu, ajnonaka@lbl.gov, 
jbbell@lbl.gov, \mbox{msday@lbl.gov}, mlminion@lbl.gov\vspace{6pt}}
\author[]{Andrew Nonaka$^{\rm \MakeLowercase{b}}$}
\author[]{John Bell$^{\rm \MakeLowercase{b}}$}
\author[]{Marcus Day$^{\rm \MakeLowercase{b}}$}
\author[]{Michael Minion$^{\rm \MakeLowercase{b}}$}
\thanks{The work at LBNL was supported by the Applied Mathematics Program
        of the DOE Office of Advanced Scientific Computing Research
        under the U.S. Department of Energy under contract DE-AC02-05CH11231.
        The submitted manuscript has been authored by a contractor of the U.S.
        Government under this contract.  Accordingly, the U.S. Government 
        retains a nonexclusive royalty-free license to publish or reproduce the 
        published form of this contribution, or allow others to do so, for U.S. 
        Government purposes. \\ \\
$^{\rm a}${\it{Division of Applied Mathematics,
              Brown University,
              Providence, RI, 02912, USA, \mbox{tel: +1 (401) 863-2115}}};\\
$^{\rm b}${\it{Center for Computational Sciences and Engineering, 
              Lawrence Berkeley National Laboratory, 
              Berkeley, CA 94720, USA, \mbox{tel: +1 (510) 486-7107}}}
}

\maketitle

\begin{abstract}
We present a fourth-order finite-volume algorithm in space and time
for low Mach number reacting flow
with detailed kinetics and transport.  Our temporal integration scheme is based on 
a multi-implicit spectral deferred correction (MISDC) strategy that iteratively couples
advection, diffusion, and reactions evolving subject to a constraint.
Our new approach overcomes a stability limitation of our previous second-order
method encountered when trying to incorporate higher-order polynomial representations 
of the solution in time to increase accuracy.  
We have developed a new iterative scheme that naturally fits within our MISDC
framework that allows us to simultaneously conserve mass and energy while satisfying
on the equation of state.  We analyse the conditions for which 
the iterative schemes are guaranteed to converge to the fixed point solution.
We present numerical examples illustrating the performance 
of the new method on premixed hydrogen, methane, and dimethyl ether flames.\\


\end{abstract}

\makeatletter{}\section{Introduction}

A broad range of problems in fluid mechanics are characterized by
dynamics with low Mach number.  In such systems, acoustic propagation
typically has negligible impact on the system state.  Low Mach number
models exploit this separation between flow and acoustic dynamics by
analytically removing sound waves from the system entirely.  In the
approximation, pressure formally becomes an elliptic field with global
coupling, and the set of conservation laws take the form of a coupled
differential-algebraic system.  Numerically, the low Mach number model
can be time-advanced on the scales of (slow) advection processes.
However, such schemes can be quite complex to implement for
multidimensional, time-dependent flows.  Moreover, there are many low
Mach number systems where the chemical and diffusive dynamics can
operate on time scales that can be much faster than the advection;
practically, the realizable time step becomes limited by the extent to
which these processes are properly coupled during time-advance.  In
this paper, we develop a highly efficient low Mach number integration
strategy that is fourth-order in space and time, while simultaneously
respecting the nonlinear coupling of all the processes.  We compare
this new algorithm to earlier versions that are both low order and
less efficiently coupled in order to demonstrate the significant
improvements afforded by the new scheme.

For smooth test problems, characterizing the convergence behavior of a
numerical integration algorithm is straighforward; with increasing
resolution higher-order spatial and temporal discretization methods
will eventually provide more accurate solutions than lower-order
methods.  For more complex problems, however, we cannot a priori
assume that the minimum resolution falls within the asymptotic range
of the methods, and therefore that higher-order methods are always
more efficient. In simple inert flow applications, we can appeal to
Reynolds numbers and Kolmogorov scales to develop absolute accuracy
requirements.  As an example of this approach, it was demonstrated in
Reference~\cite{AABM:2011} that for low-speed turbulent flows
high-order temporal and spatial discretisations outperformed
comparable second-order schemes, reducing by a factor of two in each
dimension the size of the computational mesh needed to resolve
turbulent flows at a given Reynolds number. For reacting low Mach
number flows, resolution requirements are additionally determined by
the need to accurately resolve the chemical dynamics, and to capture
the coupling between chemistry and the transport processes in the
fluid.

In \cite{LMC_SDC}  a second-order method for reacting low Mach number 
flow  with detailed kinetics and transport is introduced.
The low Mach number model is a set of differential algebraic equations
representing coupled advection, diffusion, and reaction processes that evolve subject 
to a constraint. One approach to solving such a system is to recast
the equation of state as a constraint on the
velocity divergence that determines the evolution of the thermodynamic state.
The numerical method in \cite{LMC_SDC} uses a finite-volume discretisation
in space and a time-stepping method based on
a variant of Spectral Deferred
Corrections (SDC) \cite{Dutt:2000}.
SDC is an iterative method for ordinary differential equations that has the
appealing feature that variants of arbitrarily high order can be constructed
from relatively simple lower-order methods.
In \cite{Minion:2003}, a semi-implicit variant method (SISDC) 
is introduced for ODEs with both
stiff and non-stiff processes, such as advection-diffusion systems.  The correction
equations for the non-stiff terms are discretised explicitly, whereas the stiff
corrections are treated implicitly.
Bourlioux, Layton, and Minion \cite{BLM:2003,Layton:2004} introduce a multi-implicit SDC 
approach (MISDC) for PDEs with advection, diffusion, and reaction processes. The 
advection correction equation is treated explicitly, while the diffusion and reaction 
corrections are treated implicitly and independently.  

In \cite{LMC_SDC}, a modified MISDC method is employed where 
the reaction correction equation is solved with a separate ODE solver
in a manner similar to the classic defect correction schemes \cite{zadunaisky:1964}.  
This allows the advection, diffusion, and reaction 
terms to be decoupled in the time step while retaining second-order accuracy
in time.  In addition, the tighter coupling of the terms in the deferred
corrections as compared to classical Strang splitting  results in a reduction 
of computational effort  in the reaction solves because
of smaller artificially stiff transients caused by operator splitting.
In this paper, we extend the results in \cite{LMC_SDC} to construct
a method which is fourth-order in space and time for realistic
test problems in one dimension.  The algorithm is 
closer in spirit to the original MISDC \cite{BLM:2003,Layton:2004}
approach, and  we present an analysis demonstrating why the 
approach in \cite{LMC_SDC} is ill-suited for extension to higher-order in time.  

In this paper we also develop a new technique where we evolve mass and energy
using conservation equations while satisfying the equation of state.
Our approach is based on previous
`volume discrepancy' approaches \cite{pember-flame,DayBell:2000,LMC_SDC}
where the divergence constraint is modified to drive the state variables
toward equilibrium with the ambient pressure, but now leverages the iterative 
nature of our SDC algorithm to essentially eliminate thermodynamic drift.

The rest of this paper is organised as follows.
In Section \ref{sec:equations} we review the low Mach number equation set.
In Section \ref{sec:MISDC} we review the MISDC methodology, present detailed convergence
analysis of a fourth-order variant, and discuss the implementation of this method
for a model problem.
In Section \ref{sec:methodology} we present our new volume discrepancy approach,
and describe the spatial and temporal discretisation for
the full low Mach number reacting flow equations.
In Section \ref{sec:results} we present results for our model problem as well
as several laminar flames with detailed kinetics and transport.
We summarise and conclude in Section \ref{sec:conclusions}.

\makeatletter{}\section{Low Mach Number Equation Set}\label{sec:equations}
In the low Mach number regime, the characteristic fluid velocity is small 
compared to the sound speed
(typically the Mach number is $M = U/c \sim \mathcal{O}(0.1)$ or smaller),
and the effect of acoustic wave propagation is 
unimportant to the overall dynamics of the system.  In a low Mach number 
numerical method, acoustic wave propagation is mathematically removed from the equations 
of motion, allowing for a time step based on an advective CFL condition,
\begin{equation}
\max_{\bf i} \frac{|U_i| \Delta t }{\Delta x}
\le \sigma
; \quad 0 \le \sigma \le 1,
\end{equation}
where $\sigma$ is the advective CFL number, the maximum is taken over all grid cells,
$\Delta x$ is the grid spacing, and $U_i$ is the fluid velocity in cell $i$.
Thus, this approach leads to a $\sim 1/M$
increase in the allowable time step over an explicit compressible approach.
Note that a low Mach number method does not enforce that the Mach number remain small,
but rather is suitable for flows in this regime.

In this paper, we use the low Mach number equation set from 
\cite{DayBell:2000,LMC_SDC}, which is based on the model for low Mach number combustion 
introduced by Rehm and Baum \cite{rehmBaum:1978} and rigorously derived from 
an asymptotic analysis by Majda and Sethian \cite{majdaSethian:1985}.  
We consider a gaseous mixture ignoring Soret and Dufour effects, and assume
a mixture model for species diffusion \cite{Kee:1983,Warnatz:1982}.
The resulting equations are a set of partial differential equations
representing coupled advection, diffusion, and reaction processes that are
closed by an equation of state.  The equation of state takes the form of a divergence 
constraint on the velocity, which is derived by differentiating the equation of state
in the Lagrangian frame of the moving fluid and enforcing that the
thermodynamic pressure remains constant.  Physically this manifests itself
as instantaneous acoustic equilibration to the constant thermodynamic
pressure $p_0$ (we only consider open containers in non-gravitationally stratified 
environments).  In the model, sound waves are analytically 
eliminated from our system while retaining local compressibility effects due to reactions,
mass diffusion, and thermal diffusion.

Using the notation in \cite{DayBell:2000,LMC_SDC}, the evolution equations for the 
thermodynamic variables, $(\rho,\Yb,h)$, are instantiations of mass and energy 
conservation:
\begin{eqnarray}
\frac{\partial(\rho Y_j)}{\partial t} &=& -\nabla\cdot(U\rho Y_j) + \nabla\cdot\rho\mathcal D_j\nabla Y_j + \dot\omega_j,\label{eq:cons mass}\\
\frac{\partial(\rho h)}{\partial t} &=& -\nabla\cdot(U\rho h) + \nabla\cdot\frac{\lambda}{c_p}\nabla h + \sum_j\nabla\cdot h_j\left(\rho\mathcal D_j - \frac{\lambda}{c_p}\right)\nabla Y_j,\label{eq:cons energy}
\end{eqnarray}
where $\rho$ is the density, 
$\Yb = (Y_1, \ldots, Y_N)$ are the species mass fractions,
$\mathcal D_j(\Yb,T)$ are the species mixture-averaged diffusion coefficients,
$T$ is the temperature,
$\dot\omega_j(\Yb,T)$ is the production rate for $\rho Y_j$ due to chemical reactions,
$h = \sum_j Y_j h_j$ is the enthalpy with $h_j=h_j(T)$ the enthalpy of species $j$,
$\lambda(\Yb,T)$ is the thermal conductivity, 
and $c_p = \sum_j Y_j dh_j/dT$ is the specific heat at constant pressure.
Our definition of enthalpy includes the standard enthalpy of formation, 
so there is no net change to $h$ due to reactions.
These evolution equations are closed by an equation of state,
which states that the thermodynamic pressure remain constant,
\begin{equation}
p_0 = \rho\mathcal R T\sum_j\frac{Y_j}{W_j},\label{eq:EOS}
\end{equation}
where $\mathcal R$ is the universal gas constant and $W_j$ is the molecular weight of species
$j$.  A property of multicomponent diffusive transport is that the species diffusion 
fluxes must sum to zero in order 
to conserve total mass.  For mixture models such as the 
one considered here, $\Gamma_j \equiv \rho\mathcal{D}_j\nabla Y_j$,
and that property is not satisfied in general.  Our approach is to 
identify a dominant species, in this case N$_2$, and define 
$\Gamma_{{\rm N}_2} = -\sum_{j\ne{\rm N}_2}\Gamma_j$.
Summing the species equations and noting that $\sum_jY_j = 1$ and $\sum_j\dot\omega_j=0$,
we see that (\ref{eq:cons mass}) implies the continuity equation,
\begin{equation}
\frac{\partial\rho}{\partial t} = -\nabla\cdot(U\rho).\label{eq:continuity}
\end{equation}

Equations (\ref{eq:cons mass}), (\ref{eq:cons energy}), and (\ref{eq:EOS})
form the system that we would like to solve.  Rather than directly attacking
this system of constrained differential algebraic equations, we use
a standard approach of recasting the equation of state as a divergence constraint on the
velocity field.  The constraint is derived by taking the Lagrangian derivative of equation (\ref{eq:EOS}),
enforcing that $p_0$ remain constant,
and substituting in the evolution equations for
$\rho,\Yb$, and $T$ as described in \cite{pember-flame,DayBell:2000}.
This leads to the constraint,
\begin{eqnarray}
\nabla\cdot U &=& \frac{1}{\rho c_p T}\left(\nabla\cdot\lambda\nabla T + \sum_j \Gamma_j\cdot\nabla h_j\right)\nonumber\\
&& + \frac{1}{\rho}\sum_j\frac{W}{W_j}\nabla\cdot\Gamma_j + \frac{1}{\rho}\sum_j\left(\frac{W}{W_j}-\frac{h_j}{c_p T}\right)\dot\omega_j \equiv S, \label{eq:divu}
\end{eqnarray}
where $W=(\sum_j Y_j/W_j)^{-1}$ is the mixture-averaged molecular weight.
This constraint is a linearised approximation to the velocity field required
to hold the thermodynamic pressure equal
to $p_0$ subject to local compressibility effects due to reaction heating, 
compositional changes, and thermal diffusion.

In two and three dimensions there is an evolution equation for velocity.  In second-order
schemes for incompressible and low Mach number flow, the velocity field is evolved
subject to this evolution equation, and later projected onto the vector space that
satisfies the divergence constraint.  In one dimension, the velocity field is 
uniquely specified by (\ref{eq:divu}) and the inflow boundary condition,
and therefore no projection is necessary.  We are
currently exploring ways to extend the higher-order methodologies in this paper to 
multiple dimensions subject to the divergence constraint.

\makeatletter{}\section{Multi-Implicit SDC}\label{sec:MISDC}
Here we review the SDC and MISDC methodology.
SDC methods for ODEs are introduced in Dutt et al.~\cite{Dutt:2000}.
The basic idea of SDC is to write the solution of an ODE
\begin{eqnarray}
\phi_t(t) &=& F(t,\phi(t)), \qquad t\in[t^n,t^n+\Delta t];\\
\phi(t^n) &=& \phi^n,
\end{eqnarray}
as the associated integral equation,
\begin{equation}
\phi(t) = \phi^n + \int_{t^n}^{t} F(\tau,\phi(\tau))~d\tau.
\end{equation}
We will suppress explicit dependence of $F$ and $\phi$ on $\tau$ for notational simplicity.
Given an approximation $\phi^{(k)}(t)$ to $\phi(t)$, we use an update equation to iteratively
improve the solution,
\begin{equation}
\phi^{(k+1)}(t) = \phi^n + \int_{t^n}^t \left[F(\phi^{(k+1)}) - F(\phi^{(k)})\right]d\tau +
 \int_{t^n}^t F(\phi^{(k)})~d\tau,
\end{equation}
where a low-order discretisation (e.g., forward or backward Euler) is used to approximate 
the first integral and a higher-order quadrature is used to approximate the second 
integral.  By doing so, each iteration in $k$ improves the overall order of accuracy of 
the approximation by one, up to the order of accuracy of the underlying quadrature rule 
used to evaluate the second integral.

For a given timestep, we subdivide the interval $[t^n, t^n+\Delta t]$ into $M$ subintervals,
with $M+1$ temporal nodes given by
\[
   t^n = t^{n,0} < t^{n,1} < \ldots < t^{n,M} = t^n + \Delta t \equiv t^{n+1}.
\]
For notational simplicity we will write $t^m = t^{n,m}$.
We choose the temporal nodes $t^m$ to be the appropriate Gauss-Lobatto quadrature points, 
though other choices are available \cite{Layton:2005}.  We also denote the substep time interval
by $\Delta t^m = t^{m+1} - t^m$. 
We let $\phi^{m,(k)}$ represent the $k^{\rm th}$ iterate of the solution at the 
$m^{\rm th}$ temporal node.

Bourlioux et al.~\cite{BLM:2003} and Layton and Minion \cite{Layton:2004}
introduce a variant of SDC, referred to as MISDC, in which $F$ is decomposed into distinct
processes with each treated sequentially with an appropriate explicit or 
implicit temporal discretisation.  An important difference between MISDC 
methods and operator splitting methods such as Strang splitting is that MISDC 
methods iteratively couple all physical processes together by including the effects 
of each process during the integration of any particular process.
This is in contrast to Strang splitting, 
where each process is discretised in isolation, ignoring the effects
of other processes.  Here, we write
\begin{equation}
   \phi_t = F_A(\phi) + F_D(\phi) + F_R(\phi) \equiv F(\phi),
\end{equation}
where $F_A$, $F_D$, and $F_R$ represent the advection, diffusion, and reaction 
processes respectively.  In our problems of interest, diffusion and reactions operate on 
fast time scales compared to advection.  Thus, we seek an explicit treatment of 
advection and an implicit treatment of reactions and diffusion.

We begin by initialising the solution at all temporal nodes to the solution at $t^n$,
i.e., $\phi^{m,(0)} = \phi^n,$ for all $m\in[0,M]$.
We seek to compute the next iterate of the solution, $\phi^{m,(k+1)}$ for all $m$ given that we
know $\phi^{m,(k)}$ for all $m$.  We do this by noting that $\phi^{0,(k)}=\phi^n$ for
all $k\in[0,K]$,
and then solving for each $\phi^{m+1,(k+1)}$ from $m=0$ to $M-1$ using the following sequence: 
\begin{equation}
\phi_{\rm A}^{m+1,(k+1)} = \phi^{m,(k+1)} + \int_{t^m}^{t^{m+1}}
\left[
  F_A(\phi_{\rm A}^{(k+1)}) - F_A(\phi^{(k)})
\right]dt
+ \int_{t^m}^{t^{m+1}}F(\phi^{(k)}) dt,\label{eq:MISDC_A}
\end{equation}
\begin{eqnarray}
\phi_{\rm AD}^{m+1,(k+1)} = \phi^{m,(k+1)} &+& \int_{t^m}^{t^{m+1}}
\left[
  F_A(\phi_{\rm A}^{(k+1)}) - F_A(\phi^{(k)})
+ F_D(\phi_{\rm AD}^{(k+1)}) - F_D(\phi^{(k)})
\right]dt\nonumber\\
&+& \int_{t^m}^{t^{m+1}}F(\phi^{(k)}) dt,\label{eq:MISDC_AD}
\end{eqnarray}
\begin{eqnarray}
\phi^{m+1,(k+1)} = \phi^{m,(k+1)} &+& \int_{t^m}^{t^{m+1}}
\left[
  F_A(\phi_{\rm A}^{(k+1)}) - F_A(\phi^{(k)})
+ F_D(\phi_{\rm AD}^{(k+1)}) - F_D(\phi^{(k)})\right.\nonumber\\
&&\hspace{0.5in}\left.+ F_R(\phi^{(k+1)}) - F_R(\phi^{(k)})\right]dt
+ \int_{t^m}^{t^{m+1}}F(\phi^{(k)}) dt.\nonumber\\
\label{eq:MISDC_ADR}
\end{eqnarray}
Once $\phi^{(k+1)}$ is known at all temporal nodes $m$,
the entire process can be repeated to compute the solution at all temporal nodes for the next iteration in $k$.
By using first-order discretisations in time for the first integrals in 
(\ref{eq:MISDC_A}), (\ref{eq:MISDC_AD}), and (\ref{eq:MISDC_ADR}), and
using higher-order quadrature to evaluate the second integrals, the overall accuracy of the
solution for each $k$ iterate is increased by 1, up to the order of the quadrature rule used
to evaluate the second integral over the entire time step.

In this case, if we use forward Euler to discretise advection, and backward Euler to 
discretise diffusion and reactions, we note that $\phi_{\rm A}^{m+1,(k+1)}$ does not need to 
be computed and the update consists of the following two sequential 
discretisations of equations (\ref{eq:MISDC_AD}) and (\ref{eq:MISDC_ADR}):
\begin{align}
&\begin{aligned}
\phi^{m+1,(k+1)}_{\rm AD} = \phi^{m,(k+1)} + \Delta t^m
                        \Big[&F_A(\phi^{m,(k+1)}) - F_A(\phi^{m,(k)})\\
                     + &F_D(\phi_{\rm AD}^{m+1, (k+1)})  - F_D(\phi^{m+1,(k)})\Big]
                        + I_{m}^{m+1}\left[F(\phi^{(k)})\right], \label{eq:AD_discr}
\end{aligned}\\
&\begin{aligned}
   \phi^{m+1,(k+1)} = \phi^{m,(k+1)} + \Delta t^m
                        \Big[&F_A(\phi^{m,(k+1)}) - F_A(\phi^{m,(k)})\\
                        +& F_D(\phi_{\rm AD}^{m+1, (k+1)})  - F_D(\phi^{m+1,(k)})\\
                        +& \left.F_R(\phi^{m+1, (k+1)})  - F_R(\phi^{m+1,(k)})\right] + I_{m}^{m+1}\left[F(\phi^{(k)})\right].
\label{eq:ADR_discr}
\end{aligned}
\end{align}
The second integrals in (\ref{eq:MISDC_AD}) and (\ref{eq:MISDC_ADR}) have been replaced
with numerical quadrature integrals over the substep, denoted $I_m^{m+1}$.
We note that
in \cite{BLM:2003} it was demonstrated that integration errors can be reduced by 
further subdividing the $M$ subintervals into additional nested subintervals
to treat fast-scale processes, such as diffusion and/or reactions.
Here, we choose to not further
subdivide beyond the original $M$ subintervals, so that we stay faithful to the
iterative scheme described by (\ref{eq:AD_discr}) and (\ref{eq:ADR_discr}).
Our results demonstrate we can simulate complex flames using an advective CFL
of $\sigma\sim 0.25$ without substepping diffusion or reactions.
Also, we note that 
for an evolution equation containing only advection (such as density), it is sufficient 
to use only (\ref{eq:MISDC_A}) to iteratively correct the solution, with associated discretisation,
\begin{equation}
\phi^{m+1,(k+1)} = \phi^{m,(k+1)} + \Delta t^m
 \left[F_A(\phi^{m,(k+1)}) - F_A(\phi^{m,(k)})\right]
+ I_{m}^{m+1}\left[F(\phi^{(k)})\right].\label{eq:A_discr}
\end{equation}

In \cite{LMC_SDC}, we developed a hybrid MISDC/classical deferred correction scheme
to solve the low Mach number equations.  The departure from the MISDC formulation
as originally proposed in \cite{BLM:2003} occurred when taking the time-derivative
of the reaction correction equation, given by (\ref{eq:MISDC_ADR}) in this paper
and equation (24) in \cite{LMC_SDC}, assuming that the iteratively-lagged
reaction terms cancelled, and opting to solve the ODE in equation 
(25) in \cite{LMC_SDC} instead of equation (\ref{eq:ADR_discr}) in this paper.  We then 
posited that higher-order temporal integration could be achieved 
by using higher-order polynomial representations of advection and diffusion during the 
reaction ODE step. In practice, we observed instability whenever polynomials of 
degree greater than zero were used, so we used a piecewise constant, time-centred
representation of advection and diffusion.  In Appendix \ref{app:stability} we present 
an analysis of the convergence of the numerical method in \cite{LMC_SDC}, and 
demonstrate that the generalisation to higher-order polynomial representations of
advection and diffusion results in highly unfavourable stability properties.
In light of this, the MISDC approach in this paper stays true to the MISDC
approach in \cite{BLM:2003}.  In the next Section we
analyse the new method, demonstrating convergence in the fourth-order case.

\subsection{Convergence Analysis of Fourth-Order MISDC}\label{subsec:misdc-stability}
The weakly coupled set of equations (\ref{eq:MISDC_A}), (\ref{eq:MISDC_AD}), 
and (\ref{eq:MISDC_ADR}) are chosen such that successive iterations are intended 
to converge to a fixed-point solution by sending the splitting error to zero. In order 
for the iterations to converge, we must clearly have
\begin{equation}
   \left| \phi^{m,(k+1)} - \phi^{m,(k)} \right| \to 0 \quad \text{as $k\to\infty$}.
\end{equation}
In order to study this \textit{convergence condition}, we consider the linear ODE
\begin{equation}
   \phi_t = a\phi + d\phi + r\phi \equiv F(\phi).
\end{equation}
The scalar quantities $a$, $d$, and $r$ will be proxies for our treatment
of advection, diffusion, and reaction in the full low Mach number code.
We can solve this differential equation with fourth-order accuracy, according to 
the prescription described above. We choose three Gauss-Lobatto nodes,
\begin{equation}
   t^{n, 0} = t^n, \qquad t^{n, 1} = t^n + \Delta t/2, \qquad t^{n, 2} = t^{n+1} = t^n + \Delta t.
   \label{eq:gauss-lobatto-3}
\end{equation}
As before, we denote $t^{n,m} = t^m$. The quadrature $I_m^{m+1}$ can then be 
computed by means of integrating the interpolating quadratic, obtaining the following formulas.
\begin{align}
   I_0^1(F) &= \frac{\Delta t}{24}\left(5F\left(t^0\right) + 8F\left(t^1\right) - F\left(t^2\right)\right)\label{eq:I01},\\
   I_1^2(F) &= \frac{\Delta t}{24}\left(-F\left(t^0\right) + 8F\left(t^1\right) + 5F\left(t^2\right)\right)\label{eq:I12}.
\end{align}
Given the solution at the beginning of a time step, $\phi^n$, we
compute the approximate solution at $t=t^n + \Delta t$ as follows. Initialising
$\phi^{m,(0)}=\phi^n$ for all $m$ and $\phi^{0,(k)}=\phi^n$ for all $k$,
we compute each successive iterate using equations 
(\ref{eq:AD_discr}) and (\ref{eq:ADR_discr}):
\begin{eqnarray}
\phi^{m+1,(k+1)}_{\rm AD} = \phi^{m,(k+1)} &+& \Delta t^m
                        \left[a\phi^{m,(k+1)} - a\phi^{m,(k)} + d\phi_{\rm AD}^{m+1, (k+1)}  - d\phi^{m+1,(k)}\right] \nonumber\\
                        &&+ I_{m}^{m+1}\left[F(\phi^{(k)})\right],
\end{eqnarray}
\begin{eqnarray}
   \phi^{m+1,(k+1)} = \phi^{m,(k+1)} &+& \Delta t^m
                        \left[a\phi^{m,(k+1)} - a\phi^{m,(k)}\right. + d\phi_{\rm AD}^{m+1, (k+1)} - d\phi^{m+1,(k)}\nonumber\\
                        && + \left.r\phi^{m+1, (k+1)}  - r\phi^{m+1,(k)}\right] + I_{m}^{m+1}\left[F(\phi^{(k)})\right],
\end{eqnarray}
at each temporal node $t^m$. Expanding these expressions, we can write
\begin{align}
\phi^{2,(k+1)} - \phi^{2,(k)} &= c_1\left(\phi^{1,(k)} - \phi^{1,(k-1)}\right)
                               + c_2\left(\phi^{2,(k)} - \phi^{2,(k-1)}\right),
\end{align}
where
\begin{align}
   c_1 &= \tfrac{8\Delta t(\Delta t(a^2 - d^2 + dr - r^2) + a + d + r)}
         {3(d\Delta t - 2)^2(r\Delta t - 2)^2},\label{eq:coef_1} \\
   c_2 &= \tfrac{\Delta t(-16(d+r) - a^2 \Delta t + a(8-6(d+r)\Delta t + 3dr\Delta t^2) + \Delta t(7r^2 + dr(26 - 9r\Delta t) + d^2(7 + 3r\Delta t(r\Delta t -3)))}
         {3(d-2)^2(r-2)^2}.\label{eq:coef_2}
\end{align}
Similar expressions can be derived for the difference $\phi^{1,(k+1)} - \phi^{1,(k)}$.
We can therefore conclude that
\begin{align*}
\left|\phi^{1,(k+1)} - \phi^{1,(k)}\right| &\leq \alpha \left|\phi^{1,(k)} - \phi^{1,(k-1)}\right| + \beta \left|\phi^{2,(k)} - \phi^{2,(k-1)}\right|,\\
\left|\phi^{2,(k+1)} - \phi^{2,(k)}\right| &\leq \gamma \left|\phi^{1,(k)} - \phi^{1,(k-1)}\right| + \delta \left|\phi^{2,(k)} - \phi^{2,(k-1)}\right|,
\end{align*}
where $\alpha, \beta, \gamma,$ and $\delta$ are algebraic expressions in terms 
of $a, d, r,$ and $\Delta t$. 
We see that a sufficient condition for successive iterations to converge is the 
condition $\alpha, \beta, \gamma, \delta < 1$. In other words, we define 
\[
   \theta(a, d, r, \Delta t) = \max\{\alpha, \beta, \gamma, \delta \},
\]
and require $\theta < 1$.

It is possible to compare the sets of parameters $a$, $d$, $r$, and $\Delta t$ 
that result in $\theta < 1$. For the sake of comparison, we set $a = 0$, and 
plot the regions $(d\Delta t, r\Delta t) \in \mathbb{R}^2$ such that $\theta < 1$ 
in Figure \ref{fig:stability}. 
Comparing the current method to that proposed in \cite{LMC_SDC} and analysed in 
Appendix \ref{app:stability}, we observe that the region of convergence of the 
current method encompasses a much larger range of parameters.

\begin{figure}
\centering
\subfigure[New MISDC method]{\includegraphics[width=.425\linewidth]{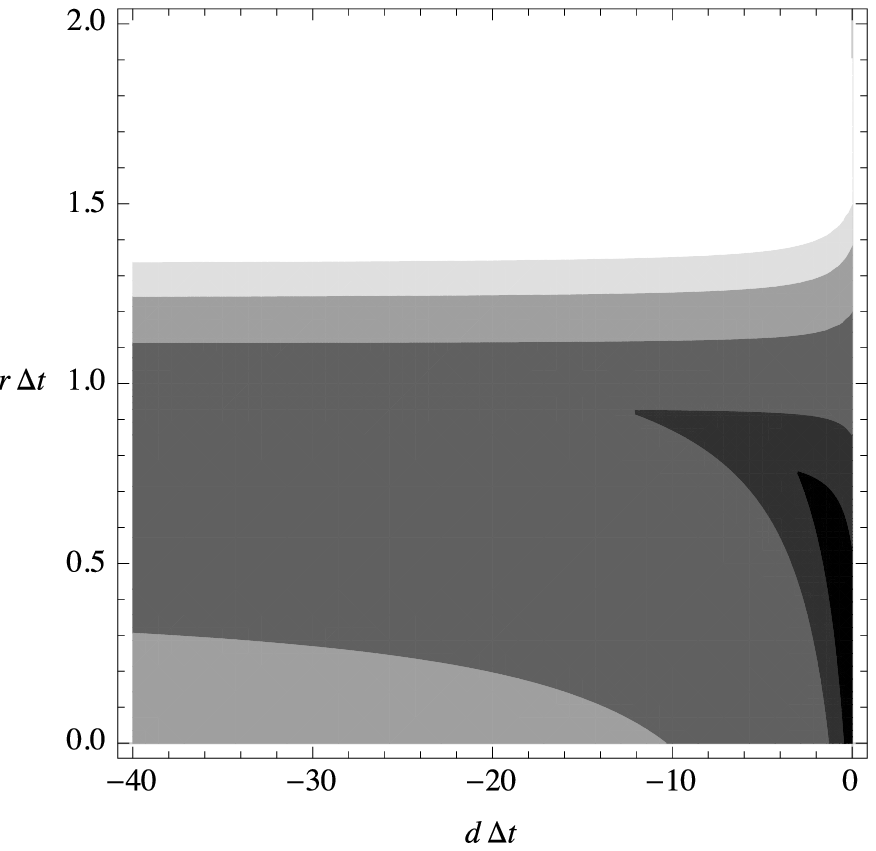}}
\subfigure[Previous proposed method]{\includegraphics[width=.565\linewidth]{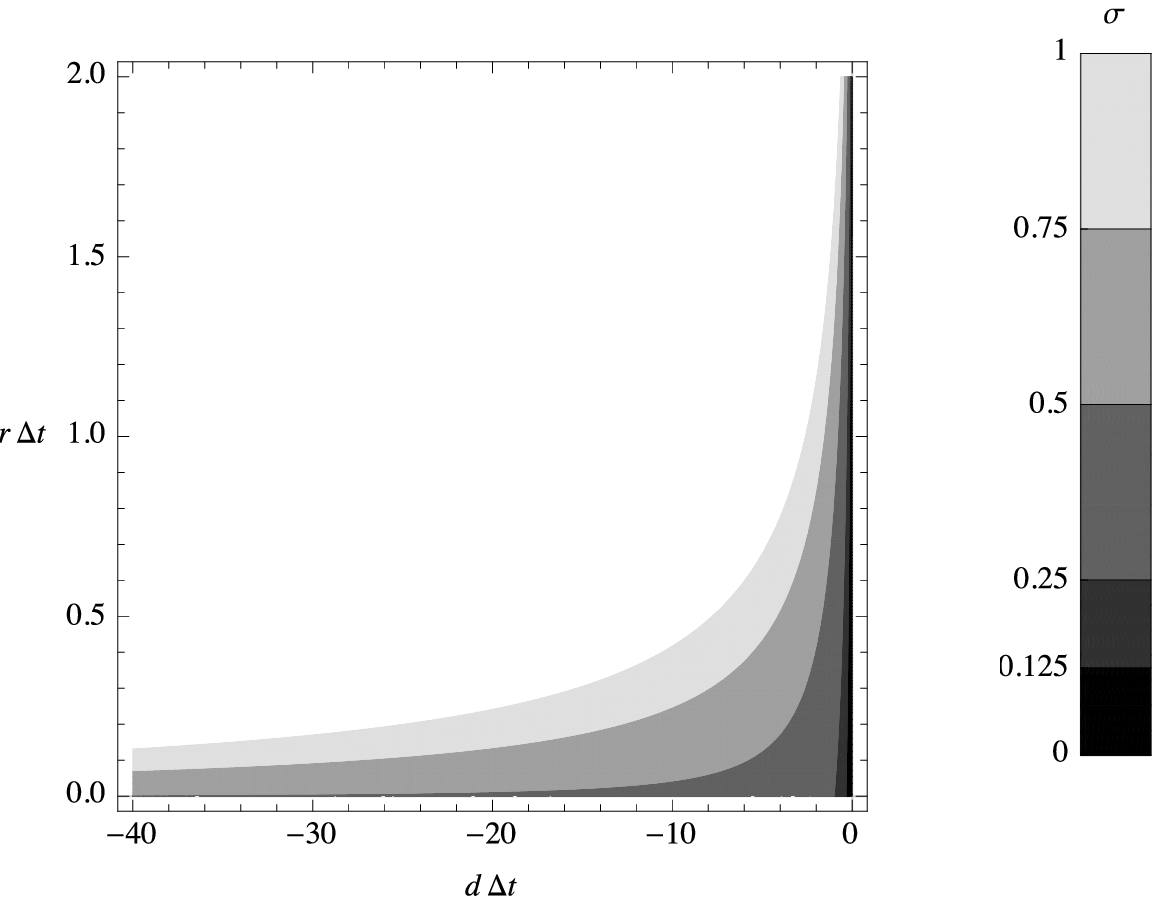}}
\caption{Convergence regions $(d \Delta t, r \Delta t) \in \mathbb{R}^2$ \label{fig:stability}}
\end{figure}

The parameters $a$, $d$ are chosen to represent the 
eigenvalues of the advection and diffusion operators, respectively. Therefore, 
$a$ can be considered to scale like $1/\Delta x$, and $d$ to scale like $1/\Delta x^2$. 
We also make the ansatz of a linear CFL, \textit{i.e.} $\Delta t = \lambda \Delta x$ for 
some $\lambda$. The reaction parameter, $r$, is independent of $\Delta x$.

We therefore write $a = \tilde a / \Delta x$, and $d = \tilde d / \Delta x^2$, 
where $\tilde a$ and $\tilde d$ are given parameters. A sufficient condition 
for the iterative scheme to converge as we send $\Delta x$ to zero is
\[
   \lim_{\Delta x \to 0} \theta(\tilde a/\Delta x, \tilde d/\Delta x^2, r, \lambda\Delta x) < 1.
\]
Calculating the limits explicitly, we see that
\begin{align*}
   \lim_{\Delta x \to 0} \alpha &= \frac{1}{12}, &    \lim_{\Delta x \to 0} \beta &= \frac{1}{3}, \\
   \lim_{\Delta x \to 0} \gamma &= \frac{2}{3}, &    \lim_{\Delta x \to 0} \delta &= \frac{7}{12},
\end{align*}
and therefore $\lim_{\Delta x \to 0} \theta = 2/3$. We can conclude that the 
fourth-order MISDC method described above is convergent in the limit as 
$\Delta x$ tends to zero. This is in contrast to the method described in 
\cite{LMC_SDC}, for which the convergence analysis is performed in 
Appendix \ref{app:stability}.  We verify fourth-order accuracy for this approach
using a test problem described in Section \ref{subsec:test-pde-results}.

\makeatletter{}\section{Numerical Methodology}\label{sec:methodology}
For the full low Mach number system we consider a one-dimensional finite volume
formulation with constant grid spacing $\Delta x$.  We describe the fourth-order
MISDC temporal integration strategy in detail in Section \ref{subsec:temporal}. 
We describe the fourth-order spatial discretisation in detail in Section \ref{spatial}.

\subsection{A New Approach for Constrained Evolution}
A major hurdle in the development of low Mach number methodologies for
complex flows is the issue that the mass and energy fields updated
with conservation equations will, in general, fail to satisfy the
equation of state.  We define $p_{\rm EOS}$ as the thermodynamic
pressure, computed directly with the equation of state using variables
updated from the conservation equations (i.e., the right-hand-side of
equation (\ref{eq:EOS})), and $p_0$ as the (constant) ambient pressure.  Our
approach of recasting the equation of state as a divergence constraint
on the velocity field is designed to constrain the evolution of the
thermodynamic state such that $p_{\rm EOS}$ remains approximately
$p_0$.  However, our formulation represents a linearisation, or
tangent-plane approximation, to the dynamics of the system.  Since the
equation of state is nonlinear, a small drift between $p_{\rm EOS}$
and $p_0$ will be observed, and indeed will grow over time.

One standard approach to resolve this issue \cite{Najm:1998,Knio:1999} is to evolve all
the thermodynamic variables but one (typically the energy or total
density field) with conservation equations, and then use the equation
of state to compute the remaining variable so that $p_{\rm EOS} = p_0$
identically.  A serious disadvantage of this approach is that it fails
to conserve mass, energy, or both.  In previous works \cite{pember-flame,DayBell:2000,LMC_SDC} we
introduced an alternative 'volume discrepancy' approach that drives
$p_{\rm EOS}$ toward $p_0$ in a way that is conservative while
maintaining the drift below a few percent.  In this paper, we exploit
the iterative nature of the MISDC advance in order to develop an
improved volume discrepancy correction.

Similar to our previous volume discrepancy approach, the constraint
equation is modified to allow additional expansion of the fluid,
accounting for the thermodynamic drift.  However, here the increment
is adjusted at each MISDC iteration, allowing us to iteratively adjust
the driving terms so that drift effectively becomes zero at the end of
each time step.  In order to construct this iteration, we return to the
derivation of the velocity constraint.  First, the equation of state,
$p=p(\rho,T,{\bf Y})$, is differentiated in the Lagrangian frame of
the moving fluid,
\begin{equation}
\frac{Dp}{Dt} = p_\rho\frac{D\rho}{Dt} + p_T\frac{DT}{Dt} + \sum_jp_{Y_j}\frac{DY_j}{Dt},\label{eq:particle paths}
\end{equation}
where the following partial derivatives are defined:
\begin{equation}
p_\rho = \left.\frac{\partial p}{\partial\rho}\right|_{T,\Yb},
\quad
p_T = \left.\frac{\partial p}{\partial T}\right|_{\rho,\Yb},
\quad
p_{Y_j} = \left.\frac{\partial p}{\partial Y_j}\right|_{\rho,T,Y_{k,k\ne j}}.
\end{equation}
Using continuity, $D\rho / Dt = -\rho \nabla \cdot U$, we rewrite (\ref{eq:particle paths}) as
\begin{equation}
\nabla\cdot U = \frac{1}{\rho p_\rho}\left(-\frac{Dp}{Dt} + p_T\frac{DT}{Dt} + \sum_jp_{Y_j}\frac{DY_j}{Dt}\right)\label{eq:particle paths2}
\end{equation}

Note that equation (\ref{eq:particle paths2}) is analytically equivalent to equation (\ref{eq:divu}) if
$Dp/Dt = \delta \chi = 0$.  Next, rather than setting $Dp/Dt = 0$, we
approximate at each cell the time derivative, $\delta \chi$, necessary
to drive the drift to zero over $\Delta t$, based on current estimates
of the advanced state:
\begin{equation}
\delta\chi = \frac{1}{\rho p_\rho}\frac{Dp}{Dt} = \frac{1}{p_0}\left(\frac{p_0 - p_{\rm EOS}}{\Delta t}\right).\label{eq:delta_chi}
\end{equation}

This field is initialised to zero on the first MISDC iteration.  After
each iteration, equation (\ref{eq:delta_chi}) is used to estimate a new correction
(increment to $\delta \chi$) required to drive the drift computed at
that iteration to zero, and this is then used in subsequent
evaluations of equation (\ref{eq:particle paths2}) for $Dp/Dt$.

The net effect of this iteration is to adjust the advection velocities
such that the conservative updates for mass and energy are both
rigorously conservative and consistent with the equation of state.  For
sufficiently resolved flows, this adjustment is small and the overall
algorithm exhibits fourth-order accuracy, as demonstrated in the
examples that follow.

\subsection{Temporal Discretisation}\label{subsec:temporal}
In our fourth-order approach, we use $M=2$ substeps 
(3 Gauss-Lobatto temporal nodes) and $K\ge 4$
MISDC iterations, but the steps below have been generalised for any number of
temporal nodes and MISDC iterations.
The steps required to advance the solution from $t^n$ to $t^{n+1}$ are as follows:\\
\begin{algorithmic}[0]
\State Set $(\rho h, \rho \Yb)^{0,(k)} = (\rho h, \rho \Yb)^n$ for all $k\in[0,K]$,
i.e., the solution at temporal node $m=0$ is a copy of the solution at $t^n$ for
all MISDC iterations.\\
\State Set $(\rho h, \rho \Yb)^{m,(0)} = (\rho h, \rho \Yb)^n$ for all $m\in[1,M]$,
i.e., the solution for MISDC iteration $k=0$ is a copy of the solution at $t^n$ for all
temporal nodes.\\
\State We use a sequence of $\delta\chi$ correction terms, one associated with 
each subinterval denoted $\delta\chi^{m-1:m,(k)}$.  Each of these terms
modifies the divergence constraint at temporal node $m$ to drive
the variables to thermodynamic equilibrium in the next MISDC iteration.
We initialise $\delta\chi^{m-1:m,(0)}=0$ for all $m\in[1,M]$.\\
\State We compute face-averaged velocities at $t^n$ by integrating the constraint
         \begin{equation}
           \nabla\cdot U^n = S^n
         \end{equation}
We set $U^{0,(k)} = U^n$ for all $k\in[0,K]$.\\
\State Next, we explicitly evaluate the right-hand-sides of the species equations
(\ref{eq:cons mass}) and enthalpy equation (\ref{eq:cons energy})
from the $k=0$ state for all $m\in[0,M]$, noting that these states are all identical
at the beginning of the time step.  These terms are used to evaluate 
$I_m^{m+1}[F(\phi^{(k)})]$ in the steps below during the first MISDC iteration.\\ \\

Now we loop over MISDC iterations (in $k$) over each temporal node (in $m$):\\
\For{$k=0 \text{ \bf to } K-1$}
   \For{$m=0 \text{ \bf to } M-1$}\\
      \State
         \parbox[t]{\dimexpr\linewidth-\algorithmicindent}
         {If $m>0$,
         correct to the divergence constraint to account for the fact that
         the $p_{\rm EOS}$ does not match $p_0$ in the most recent solution,}
         \begin{equation}
            \delta\chi^{m-1:m,(k+1)} = \delta\chi^{m-1:m,(k)} + \frac{2}{p_0}\left(\frac{p_{\rm EOS}^{m,(k+1)}-p_0}{\Delta t^{m-1}}\right),\label{eq:delta_chi update}
         \end{equation}
      \State
         \parbox[t]{\dimexpr\linewidth-\algorithmicindent}
         {where, again, $p_{\rm EOS}$ is a function of $\rho,\Yb,T$ given by the 
           right-hand-side of equation (\ref{eq:EOS}).
           Note, the factor of 2 in (\ref{eq:delta_chi update}) reflects that our correction is applied
           piecewise linearly over the time interval from $t^m$ to
           $t^{m+1}$, with zero increment of $\delta \chi$ at $t^m$.
           We then compute face-averaged velocities, $U^{m,(k+1)}$, by integrating the 
           constraint,}
         \begin{equation}
           \nabla\cdot U^{m,(k+1)} = S^{m,(k+1)} + \delta\chi^{m-1:m,(k+1)}.
         \end{equation}
      \State
         \parbox[t]{\dimexpr\linewidth-\algorithmicindent}
         {(For all $m$) Compute $\rho^{m+1,(k+1)}$ by discretising the continuity equation 
         (\ref{eq:continuity}) with the MISDC correction equation for PDEs containing
         only advection terms (\ref{eq:A_discr}),}
         \begin{eqnarray}
            \rho^{m+1,(k+1)} = \rho^{m,(k+1)} &+& \Delta t^m\left[-\nabla\cdot(U\rho)^{m,(k+1)}
            + \nabla\cdot(U\rho)^{m,(k)}\right] \nonumber\\
             &+& I_m^{m+1} \left[-\nabla\cdot(U\rho)^{(k)}\right].
             \label{eq:rho-update}
         \end{eqnarray}
      \State
         \parbox[t]{\dimexpr\linewidth-\algorithmicindent}
         {Compute updated mass fractions, $Y_{j,{\rm AD}}^{m+1,(k+1)}$,
         by discretising the species equations (\ref{eq:cons mass}) with the MISDC
         advection-diffusion correction equation (\ref{eq:AD_discr}).
         This amounts to solving the implicit system,}
         \begin{alignat}{3}
         \rho^{m+1,(k+1)} Y_{j,{\rm AD}}^{m+1,(k+1)} &= \mathmakebox[2em][l]{(\rho Y_j)^{m,(k+1)}} \nonumber\\
         & + \Delta t^m \Big[&&-\nabla\cdot(U\rho Y_j)^{m,(k+1)} + \nabla\cdot(U\rho Y_j)^{m,(k)}\nonumber\\
         &&& + \nabla\cdot\rho^{m+1,(k)}\mathcal{D}_j^{m+1,(k)}\nabla Y_{j,{\rm AD}}^{m+1,(k+1)} - \nabla\cdot\Gamma_j^{m+1,(k)}\Big]\nonumber\\
         & \mathmakebox[2em][l]{+ I_m^{m+1}\left[-\nabla\cdot(U\rho Y_j) + \nabla\cdot\Gamma_j + \dot\omega\right]^{(k)}.}
         \end{alignat}
         \State Define $\Gamma_{j,{\rm AD}}^{m+1,(k+1)} = \rho^{m+1,(k)}\mathcal{D}_j^{m+1,(k)}\nabla Y_{j,{\rm AD}}^{m+1,(k+1)}$.
      \State
         \parbox[t]{\dimexpr\linewidth-\algorithmicindent}
         {Compute updated enthalpy, $h_{\rm AD}^{m+1,(k+1)}$, by discretising the enthalpy
         equation (\ref{eq:cons energy}) with the MISDC
         advection-diffusion correction equation (\ref{eq:AD_discr}).  
         We remark that the differential diffusion terms are treated explicitly in 
         order to avoid a more complicated linear system.
         This amounts to solving the implicit system,}
         \begin{alignat}{3}
            &\hspace{-0.5in}\mathmakebox[2em][l]{\rho^{m+1,(k+1)}h_{\rm AD}^{m+1,(k+1)} = (\rho h)^{m,(k+1)}} \nonumber\\
            & + \Delta t^m \Bigg[&&-\nabla\cdot(U\rho h)^{m,(k+1)} + \nabla\cdot(U\rho h)^{m,(k)} \nonumber\\
            &&& + \nabla\cdot\frac{\lambda^{m+1,(k)}}{c_p^{m+1,(k)}}\nabla h_{\rm AD}^{m+1,(k+1)} - \nabla\cdot\frac{\lambda^{m+1,(k)}}{c_p^{m+1,(k)}}\nabla h^{m+1,(k)} \nonumber\\
            &&& + \sum_j\nabla\cdot h_j^{m,(k+1)}\left(\Gamma_j^{m,(k+1)} - \frac{\lambda^{m,(k+1)}}{c_p^{m,(k+1)}}\nabla Y_j^{m,(k+1)}\right) \nonumber\\
            &&& - \sum_j\nabla\cdot h_j^{m,(k)}\left(\Gamma_j^{m,(k)} - \frac{\lambda^{m,(k)}}{c_p^{m,(k)}}\nabla Y_j^{m,(k)}\right) \Bigg] \nonumber\\
            & \mathmakebox[2em][l]{ + I_m^{m+1}\left[-\nabla\cdot(U\rho h) + \nabla\cdot\frac{\lambda}{c_p}\nabla h + \sum_j\nabla\cdot h_j\left(\Gamma_j - \frac{\lambda}{c_p}\nabla Y_j\right)\right]^{(k)}.}
         \end{alignat}
      \State
      \State
         \parbox[t]{\dimexpr\linewidth-\algorithmicindent}{
         Note that since there is no contribution due to reactions in the 
         enthalpy update, we can say}
         \begin{equation}
            (\rho h)^{m+1,(k+1)} = \rho^{m+1,(k+1)}h_{\rm AD}^{m+1,(k+1)}.
            \label{eq:rhoh-reac}
         \end{equation}
      \State
         \parbox[t]{\dimexpr\linewidth-\algorithmicindent}
         {Next, we solve the reaction correction equation for $(\rho Y_j)^{m+1,(k+1)}$
         using the MISDC advection-diffusion-reaction correction equation
         (\ref{eq:ADR_discr}),}
         \begin{eqnarray}
         (\rho Y_{j})^{m+1,(k+1)} &=&\mathmakebox[2em][l]{(\rho Y_j)^{m,(k+1)}} \nonumber\\
         &&+ \Delta t^m \Big[ -\nabla\cdot(U\rho Y_j)^{m,(k+1)} - \nabla\cdot(U\rho Y_j)^{m,(k)}\nonumber\\
         &&\hspace{0.5in} +\nabla\cdot\Gamma_{j,{\rm AD}}^{m+1,(k+1)} - \nabla\cdot\Gamma_j^{m+1,(k)}\nonumber\\
         &&\hspace{0.5in} +\dot\omega^{m+1,(k+1)} - \dot\omega^{m+1,(k)}\Big]\nonumber\\
         && \mathmakebox[2em][l]{ + I_m^{m+1}\left[-\nabla\cdot(U\rho Y_j) + \nabla\cdot\Gamma_j + \dot\omega\right]^{(k)}. }
         \label{eq:rhoY-reac}
         \end{eqnarray}
      \State \parbox[t]{\dimexpr\linewidth-\algorithmicindent}{
         See Section \ref{Sec:Reactions} for our solution technique for these nonlinear implicit equations.\\}
      \State
   \State \EndFor~(end loop over temporal nodes $m$)\\
   \State \parbox[t]{\dimexpr\linewidth-\algorithmicindent}
     {Since in general, $p_{\rm EOS}$ at the final temporal node (the ``$M,(k+1)$" 
       state) is not in thermodynamic 
       equilibrium with $p_0$, we correct the divergence constraint.
       We do this by incrementing $\delta\chi^{M-1:M,(k+1)}$ using}
   \begin{equation}
     \delta\chi^{M-1:M,(k+1)} = \delta\chi^{M-1:M,(k)} + \frac{2}{p_0}\left(\frac{p_{\rm EOS}^{M,(k+1)}-p_0}{\Delta t^{M-1}}\right).
   \end{equation}
   \State \parbox[t]{\dimexpr\linewidth-\algorithmicindent}{
          Compute $U^{M,(k+1)}$ by integrating the constraint, 
          \begin{equation}
            \nabla\cdot U^{M,(k+1)} = S^{M,(k+1)} + \delta\chi^{M-1:M,(k+1)},
          \end{equation}
          and evaluate the right-hand-sides of the species equations
          (\ref{eq:cons mass}) and enthalpy equation (\ref{eq:cons energy})
          from the '$M,(k+1)$' state.\\}
\EndFor~(end loop over MISDC iterates $k$) \\
\State Advance the solution by setting $(\rho h, \rho \Yb)^{n+1} = (\rho h, \rho \Yb)^{M,(K)}$.\\
\end{algorithmic}
\subsubsection{Solving the Reaction Correction Equations}\label{Sec:Reactions}
To solve the reaction correction equations (\ref{eq:rhoY-reac}) for $(\rho Y_m)^{m+1,(k+1)}$,
we use Newton's method for this implicit system.  Note that since there 
is no reaction contribution to the enthalpy equation, $(\rho h)^{m+1,(k+1)}$
is already known from (\ref{eq:rhoh-reac}). Likewise, the density $\rho^{m+1,(k+1)}$ is 
known from equation (\ref{eq:rho-update}).  The production rates are
a function of $\Yb$ and $h$.  Writing out the mass fractions
$\Yb \equiv (Y_1, \ldots, Y_N)$,
and the production rates $\dot{\bm\omega} \equiv (\dot\omega_1, \ldots, \dot\omega_N)$,
equation (\ref{eq:rhoY-reac}) takes the form of a nonlinear backward 
Euler-type equation for $\Yb$:
\begin{equation}
   \label{eq:be-reac}
   \rho^{m+1,(k+1)}\Yb - \Delta t^m \dot{\bm\omega}(\Yb) = {\bf b}.
\end{equation}
We use an analytic Jacobian \cite{Perini:2012}
and using the solution from the previous MISDC iterate as an initial guess.
The Newton solve has been observed to converge to within tolerance within only a 
few iterations for all $k$, and in subsequent MISDC iterations, the initial guess
improves with each iteration and even fewer Newton iterations are needed.
We iterate until the max norm of the residual is less than $10^{-14}$.  In our testing,
using an even tighter tolerance of $10^{-16}$ did not affect the convergence rates
in the significant figures we report below.

\subsection{Spatial Discretisation}\label{spatial}
We use a finite volume discretisation with $n$ cells, indexed from $i=(0,\ldots,n-1)$.
We distinguish between three types of quantities.  A cell-averaged quantity is denoted
by angle brackets,
\[
\langle\phi\rangle_i \equiv
\frac{1}{\Delta x} \int_{x_{i-\half}}^{x_{i+\half}} \phi(x)~dx.
\]
A cell-centred quantity is
denoted by a hat, $\widehat{\phi}_i \equiv \phi(x_i)  $, and a face-averaged quantity is
denoted by a tilde, $\widetilde\phi_{i+\half} \equiv \phi(x_{i+\half})$.  Note that in one 
dimension, face-averages are simply point values at the endpoints of a cell.

To convert between these values, as well as compute fourth-order
gradients, products, and quotients, we rely on standard operations found in the
finite volume literature \cite{AABM:2011,Kadioglu:2008,McCorquodale:2011,Zhang:2012}.
The fourth-order formulas to convert from cell-averaged to 
cell-centred, and vice versa, are,
\begin{eqnarray}
\label{eq:avg-to-center}
\widehat\phi_i &=& \langle\phi\rangle_i - \frac{1}{24}(\langle\phi\rangle_{i-1}-2\langle\phi\rangle_i+\langle\phi\rangle_{i+1}),\\
\label{eq:center-to-avg}
\langle\phi\rangle_i &=& \widehat\phi_i + \frac{1}{24}(\widehat\phi_{i-1}-2\widehat\phi_i+\widehat\phi_{i+1}).
\end{eqnarray}
We can compute a fourth-order approximation of a quantity at cell faces given either 
cell-centred values or cell-averaged values using the following stencils:
\begin{align}
\label{eq:avg-to-face}
\widetilde\phi_{i+\half} &= \frac{-\langle\phi\rangle_{i-1} + 7\langle\phi\rangle_i + 7\langle\phi\rangle_{i+1} - \langle\phi\rangle_{i+2}}{12},\\
\label{eq:center-to-face}
\widetilde\phi_{i+\half} &= \frac{-\widehat\phi_{i-1} + 9\widehat\phi_i + 9\widehat\phi_{i+1} - \widehat\phi_{i+2}}{16}.
\end{align}
The fourth-order approximation to the gradient at a cell face is
\begin{equation}
\label{eq:grad-from-avg}
\widetilde{\nabla\phi}_{i+\half} = \frac{\langle\phi\rangle_{i-1} - 15\langle\phi\rangle_i + 15\langle\phi\rangle_{i+1} - \langle\phi\rangle_{i+2}}{12\Delta x}.
\end{equation}
Given cell-averaged quantities $\langle \phi \rangle_i$ and $\langle \psi \rangle_i$, we can compute a fourth-order 
approximation to the cell-average of the product by
\begin{equation}
   \label{eq:avg-product}
   \langle \phi \psi \rangle_i = \langle \phi \rangle_i \langle \psi \rangle_i + \frac{\Delta x^2}{12} \phi^{\rm G}_i \psi^{\rm G}_i + \mathcal{O}(\Delta x^4),
\end{equation}
where $\phi^{\rm G}_i$ and $\psi^{\rm G}_i$ are given by the gradient formula, e.g.,
\begin{equation}
   \phi^{\rm G}_i = \frac{5 \langle\phi\rangle_{i-2} - 34 \langle\phi\rangle_{i-1} + 34 \langle\phi\rangle_{i+1} - 5 \langle\phi\rangle_{i+2}}{48 \Delta x}.
\end{equation}
Similarly, we can compute a fourth-order approximation to the cell-average of a quotient by
\begin{equation}
  \label{eq:avg-quotient}
  \left\langle\frac{\phi}{\psi}\right\rangle_i = \frac{\langle\phi\rangle_i}{\langle\psi\rangle_i}
     + \frac{\Delta x^2}{12}\left(
     \frac{\langle\phi\rangle_i \left(\psi^{\rm G}_i\right)^2}{\langle\psi\rangle_i^3} - 
     \frac{\phi^{\rm G}_i \psi^{\rm G}_i}{\langle\psi\rangle_j^2}
     \right) + \mathcal{O}(\Delta x^4).
\end{equation}
At inflow and outflow, the strategy is to use the boundary condition and four interior
data values to extrapolate two ghost cells values to fourth-order accuracy, allowing
us to use these same stencils.
At inflow, we have the Dirichlet value at the face, $\phi_b$.  Given cell-averaged data,
the ghost cell-averaged values are
\begin{equation}
\langle\phi\rangle_{-1} = \frac{60\phi_b - 77\langle\phi\rangle_0 + 43\langle\phi\rangle_1 - 17\langle\phi\rangle_2 + 3\langle\phi\rangle_3}{12}
\label{eq:ghost-avg-1}
\end{equation}
\begin{equation}
\langle\phi\rangle_{-2} = \frac{300\phi_b - 505\langle\phi\rangle_0 + 335\langle\phi\rangle_1 - 145\langle\phi\rangle_2 + 27\langle\phi\rangle_3
}{12}
\label{eq:ghost-avg-2}
\end{equation}
Given cell-centred data, the ghost cell-centred values are
\begin{equation}
\widehat{\phi}_{-1} = \frac{128\phi_b - 140\widehat{\phi}_0 + 70\widehat{\phi}_1 - 28\widehat{\phi}_2 + 5\widehat{\phi}_3}{35}
\end{equation}
\begin{equation}
\widehat{\phi}_{-2} = \frac{128\phi_b - 210\widehat{\phi}_0 + 140\widehat{\phi}_1 - 63\widehat{\phi}_2 + 12\widehat{\phi}_3}{7}
\end{equation}

At outflow, we have a homogeneous Neumann condition.   Given cell-averaged data,
the ghost cell-averaged values are
\begin{equation}
\langle\phi\rangle_n = \frac{5\langle\phi\rangle_{n-1} + 9\langle\phi\rangle_{n-2} - 5\langle\phi\rangle_{n-3} + \langle\phi\rangle_{n-4}}{10}
\label{eq:ghost-avg n}
\end{equation}
\begin{equation}
\langle\phi\rangle_{n+1} = \frac{-15\langle\phi\rangle_{n-1} + 29\langle\phi\rangle_{n-2} - 15\langle\phi\rangle_{n-3} + 3\langle\phi\rangle_{n-4}}{2}
\label{eq:ghost-avg n+1}
\end{equation}
Given cell-centred data, the ghost cell-centred values are
\begin{equation}
\widehat{\phi}_n = \frac{17\widehat{\phi}_{n-1} + 9\widehat{\phi}_{n-2} - 5\widehat{\phi}_{n-3} + \widehat{\phi}_{n-4}}{22}
\end{equation}
\begin{equation}
\widehat{\phi}_{n+1} = \frac{-135\widehat{\phi}_{n-1} + 265\widehat{\phi}_{n-2} - 135\widehat{\phi}_{n-3} + 27\widehat{\phi}_{n-4}}{22}
\end{equation}
To compute the advection terms, we use the divergence theorem,
\begin{equation}
\langle \nabla\cdot(U\phi)\rangle_i = \frac{\widetilde U_{i+\half}\widetilde\phi_{i+\half} - \widetilde U_{i-\half}\widetilde\phi_{i-\half}}{\Delta x},
\end{equation}
where $\widetilde\phi$ on faces is computed using formula (\ref{eq:avg-to-face}).
We use the inflow boundary condition and
integrate the divergence constraint to obtain the velocity,
\begin{equation}
\widetilde U_{i+\half} - \widetilde U_{i-\half} = \langle S + \delta\chi\rangle_i\Delta x .
\end{equation}
Note that the terms comprising $S + \delta\chi$ are initially computed
at cell centres, and then
converted to a cell-average using formula (\ref{eq:avg-to-center}).

The diffusion operators from equations (\ref{eq:cons mass}) and 
(\ref{eq:cons energy}) are $\nabla\cdot \rho\mathcal{D}_j\nabla Y_j$ and 
$\nabla\cdot (\lambda/c_p)\nabla h$. The implicit linear solve for these 
operators takes the general form
\[
 \rho \phi - \nabla\cdot  D \nabla \phi = b,
\]
where $ D$ represents the diffusion coefficient. The finite volume 
discretisation of this solve can be written
\begin{equation}
   \langle \rho \phi \rangle - \left\langle \nabla\cdot D\nabla \phi \right\rangle = \langle b \rangle.\label{eq:diffusion_discr}
\end{equation}
Equation (\ref{eq:diffusion_discr}) can be framed as a linear solve for the 
cell-average $\langle \phi \rangle$ as follows.
The first term on the left-hand side is the cell-average of the product of 
$\rho$ and $\phi$, which we compute using the product rule (\ref{eq:avg-product}). 
The second term is computed as
\begin{equation}
	\left\langle \nabla\cdot D\nabla \phi \right\rangle =
		\frac{\widetilde{ D}_{i+\half} \widetilde{\nabla\phi}_{i+\half} - \widetilde{ D}_{i-\half} \widetilde{\nabla\phi}_{i-\half}}
		{\Delta x}.
\end{equation}
The gradient of $\phi$ on faces is computed using formula (\ref{eq:grad-from-avg}).
The diffusion coefficient $ D$ must be computed at faces.
We use cell-centred $\widehat \rho$, $\widehat{\rho h}$, and 
$\widehat{\rho Y_j}$, in order to compute cell-centred diffusion coefficients
$\widehat{\mathcal{D}_j}$, $\widehat \lambda$, and $\widehat {c_p}$.
These cell-centred values are then averaged to faces using (\ref{eq:center-to-face}). 
Near the boundary, the entries of the matrix must be modified to take into account 
the ghost cells, whose values are defined by equations 
(\ref{eq:ghost-avg-1})-(\ref{eq:ghost-avg-2}) and 
(\ref{eq:ghost-avg n})-(\ref{eq:ghost-avg n+1}). At the inflow boundary, since the
inflow value is inhomogeneous, the right-hand-side, $\langle b \rangle$ 
must also be modified to respect the Dirichlet condition.  The resulting matrix is banded, 
and is pentadiagonal, except in the first and last rows, which include an additional term.
After solving the banded linear system for $\langle \phi \rangle$, we use the product rule
to compute the solution $\langle \rho \phi \rangle.$

Concerning the the backward-Euler equation for reactions (\ref{eq:rhoY-reac}),
we convert the right-hand-side to cell-centred values using equation (\ref{eq:avg-to-center}). 
We then perform the nonlinear backward Euler solve detailed in Section \ref{Sec:Reactions}. 
From this solution, we can then obtain the cell-centred production rates, $\widehat{\dot\omega_j}$. 
The production rates are then converted to cell-averaged values, $\langle \dot\omega_j \rangle$ 
using equation (\ref{eq:center-to-avg}). This term is then substituted into equation (\ref{eq:rhoY-reac}) 
in order to update $\langle \rho Y_j \rangle$ using a cell-averaged right-hand-side.

\makeatletter{}\section{Results}\label{sec:results}
In this section, we present results both for a test problem
using the algorithm in Section \ref{subsec:misdc-stability}
and for the one-dimensional low Mach number combustion
algorithm, simulating three types of premixed laminar flames
with detailed kinetics and transport (hydrogen, methane, and dimethyl ether).
We verify fourth-order accuracy in all of these cases.

In the following tests, we perform the simulations at various resolutions,
decreasing $\Delta x$ by a factor of two, while holding the advective CFL
number constant. We estimate the error by comparing the solution at resolution
$\Delta x$ with the solution at resolution $\Delta x / 2$. For a simulation at
a coarse resolution with $n_c$ cells, we compute the $L^1$ error using
\begin{equation}
   L^1_{n_c} = \frac{1}{n_c}\sum_{i=1}^{n_c} \left| \phi_i^c - \phi_i^{f\rightarrow c} \right|,
\end{equation}
where $\phi^c$ is the coarse solution,
$\phi^{f\rightarrow c}$ is a coarsened version of the solution with twice the
resolution ($n_f=2n_c$ cells).  For the finite-difference test PDE, the coarsening
is done by direct injection, and for our finite volume flame simulations, we
average the fine solution to the coarser grid.
We can then define the convergence rate by
\begin{equation}
   r^{n_c/n_f} = \log_2\left(\frac{L^1_{n_c}}{L^1_{n_f}}\right).
\end{equation}

\subsection{Test PDE}\label{subsec:test-pde-results}
As a test bed for the MISDC method described in Section \ref{subsec:misdc-stability},
we consider the initial boundary value problem
\begin{equation}
   \label{eq:test-pde}
   \left\{
   \begin{alignedat}{5}
      \phi_t(x,t) &= a\phi_x + d \phi_{xx} + r\phi(\phi-1)(\phi-\myhalf) \qquad\text{for } (x,t) \in [0,20]\times[0,T],\\
      \phi(0,t) & = 1,\\
      \phi(20,t) &= 0,\\
      \phi(x,0) &= \phi^0(x),
   \end{alignedat}\right.
\end{equation}
We choose the initial condition to be given by
\begin{equation}
   \phi^0(x) = \frac{\tanh(10-2x)+1}{2}.
\end{equation}
We can then solve this equation using the the method of 
lines. The advection term $a\phi_x$ is approximated by a fourth-order finite difference 
operator $A(\phi)$. The diffusion term $d \phi_{xx}$ is approximated by a 
fourth-order Laplacian operator, denoted $D(\phi)$. Both operators are chosen 
to respect the Dirichlet boundary conditions.
We denote $R(\phi) = r\phi(\phi-1)(\phi-\myhalf)$. Thus, our PDE has the form,
\begin{equation}
   \phi_t = A(\phi) + D(\phi) + R(\phi),
\end{equation}
which is solved using the MISDC method described in Section \ref{subsec:misdc-stability}

We begin by setting the solution for the $k=0$ iterate at all temporal nodes to
the solution at $t^n$,
i.e., $\phi^{m,(0)} = \phi^n$ for all $m\in[0,M]$.
We treat the advective process explicitly, and the diffusion and reaction 
processes implicitly.  As noted in \cite{BLM:2003}, since advection is treated 
explicit, we do not need to compute a provisional solution for advection, $\phi_A^{m,(k+1)}$.
At each temporal node $t^{m}$, we compute a provisional 
solution $\phi^{m,(k+1)}_{AD}$ by performing a sparse, banded linear solve.
This provisional solution is then used in the correction equation for 
the updated solution, $\phi^{m,(k+1)}$.  In order to solve the reaction correction equation, 
we use Newton's method.  The solution from the previous iterate is chosen to be the initial 
guess for the Newton solver.

\begin{figure}
\begin{center}
\includegraphics[width=5in]{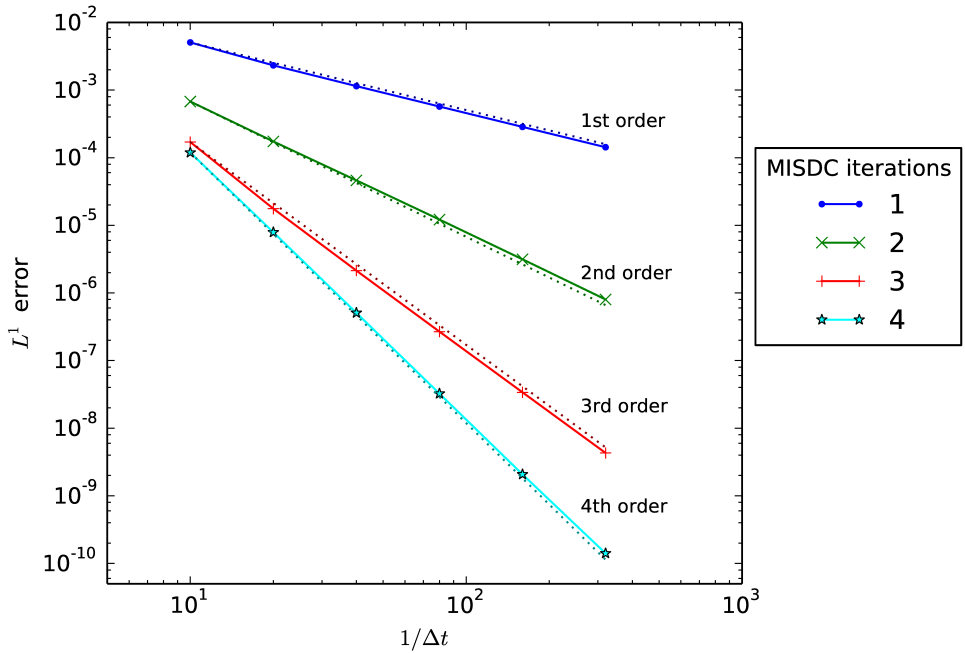}
\caption{\label{fig:test-pde-results} Log-log plot of $L^1$ error versus 
$1/\Delta t$ for the test PDE (\ref{eq:test-pde}) using three Gauss-Lobatto 
nodes. Solid lines are the numerical results, and dotted lines indicate the 
expected convergence results.}
\end{center}
\end{figure}

Using the methods from Section \ref{subsec:misdc-stability}, we obtain the expected order of 
accuracy, given by $\min\{K, Q\}$, where $K$ is the number of MISDC iterations, 
and $Q$ is the order of quadrature. Using three Gauss-Lobatto nodes, the 
quadrature is fourth-order accurate. Therefore, the overall order of accuracy is 
equal to the number of MISDC iterations, up to a maximum of four. We set the 
parameters $a = -0.1$, $d = 1$, and $r = -10$. We start with an initial 
discretisation of $n=200$ gridpoints, and set $\Delta t = \Delta x / 2$. 
Simultaneously refining in space and time we obtain the results for $L^1$ error 
shown in Figure \ref{fig:test-pde-results}.  Each numerical test
gives the expected order of accuracy, up to fourth-order

\subsection{Flame Simulations}
We now analyse the convergence behaviour of our scheme on a set of more complex 
problems: premixed laminar flames burning hydrogen, methane, and dimethyl ether
fuels that propagate through the domain.  In all cases here, the diffusive and 
reactive processes are well-known to
be numerically 'stiff' on the advection time scales used to set our numerical time step.
Because of this, we find that eight MISDC iterations are required per time step for robust 
integration, and we demonstrate in the following section that our strategy with this setting
results in the expected convergence properties.

\subsubsection{Hydrogen Flame}\label{sec:hydrogen flame}

We study the performance of the MISDC algorithm to propagate
a one-dimensional premixed hydrogen flame.
The simulations are based on the GRIMech-3.0 \cite{GRI30:1995} 
model and associated databases for thermodynamic relationships and mixture-averaged
transport properties, as given in the CHEMKIN-III library \cite{CHEMKIN:1996}
format. Note that we manually stripped the carbon-containing species and
associated reactions from the model, since they are irrelevant for the hydrogen-air case.
The resulting model for hydrogen-air mixtures consists of 9 species and 27 reactions.  

The detailed structure of these flames feature the prominent role of molecular and atomic 
hydrogen, both of which diffuse considerably faster than the other species in the
system.  This differential diffusion between species
has significant impact on steady 1D profiles; H and H$_2$ profiles tend to
be considerably broader than the others, and this plays a key role in the 
flame stabilisation.

In this configuration, an unstrained one-dimensional flame propagates into a homogeneous
hydrogen-air mixture.  A steady solution consists of thermal
and species profiles co-moving in a frame with the flame
propagation.  In the frame of the unburned fuel, the steady solution
propagates toward the inlet at the unstrained laminar burning speed, $s_L$,
which is a function of the inlet state.  At the chosen conditions, 
$Y({\rm H}_2 : {\rm O}_2 : {\rm N}_2) = (0.0107 : 0.2304 : 0.7589)$, $p$=1~atm,
and $T$=298~K, $s_L = 14.9$~cm~s$^{-1}$.  For each of our fourth-order
flame simulations, we arbitrarily set the inlet velocity to 5~cm/s.
Thus, the hydrogen flame propagates toward the inlet at 9.9~cm~s$^{-1}$.

Initial flame profiles are generated for this study in two auxiliary steps.
First, a steady one-dimensional solution is computed using the
PREMIX code \cite{PREMIX:1985}.  PREMIX incorporates a first-order difference 
scheme on nonuniform grid in one dimension.  The PREMIX solution is translated
into the frame of the unburned fuel, and interpolated onto a uniform grid
with 8192 cells across a 1.2~cm domain.  While this solution exhibits the 
essential features of the flame, it is not C2-continuous; higher-order discontinuities
will pollute subsequent convergence analysis.  To resolve this issue, we use the 
second-order low Mach number code from our previous work \cite{LMC_SDC},
to evolve the PREMIX solution for an additional $160~\mu$s using a time step of
$\Delta t = 0.8~\mu$s.  Finally, the initial data to test our fourth-order 
algorithm is generated by averaging this solution to a set of uniform meshes,
using $n =$ 128, 256, 512, and 1024 cells.  

The initial data at each resolution are evolved for 1.6~ms so that the flame propagates
approximately 16 $\mu$m across the mesh at a $\Delta t$ corresponding to an advective
CFL of $\sigma\approx$ 0.28.  The resulting profiles are compared to a reference solution
as discussed above. The error and convergence results are presented in Table \ref{table:hydrogen}.
Fourth-order convergence is observed in all variables.  Note that HO$_2$ and H$_2$O$_2$ have the narrowest
profiles of all species, and are therefore the most demanding to converge.

\begin{table}[tb]
\caption{Error and convergence rates for a premixed hydrogen flame using 
the fourth-order MISDC method with an advective CFL of $\sigma\approx 0.28$.
	\label{table:hydrogen}}
\vspace{10pt}
\begin{center}
	\begin{tabular}{ l c c c c c }
	\hline
	Variable & $L^1_{128}$ & $r^{128/256}$ & $L^1_{256}$ & $r^{256/512}$ & $L^1_{512}$\\
	\hline\\[-1.5ex]
$Y({\rm H}_2)$          & 5.91E-08 & 4.01 & 3.67E-09 & 3.98 & 2.33E-10 \\
$Y({\rm O}_2)$          & 1.10E-06 & 4.00 & 6.83E-08 & 4.05 & 4.14E-09 \\
$Y({\rm H}_2{\rm O})$   & 1.01E-06 & 4.01 & 6.25E-08 & 4.05 & 3.76E-09 \\
$Y({\rm H})$            & 1.17E-09 & 3.70 & 9.00E-11 & 3.91 & 5.97E-12 \\
$Y({\rm O})$            & 2.70E-08 & 3.93 & 1.77E-09 & 4.01 & 1.10E-10 \\
$Y({\rm OH})$           & 3.17E-08 & 4.01 & 1.97E-09 & 4.06 & 1.18E-10 \\
$Y({\rm HO}_2)$         & 3.56E-08 & 3.71 & 2.72E-09 & 3.88 & 1.86E-10 \\
$Y({\rm H}_2{\rm O}_2)$ & 1.41E-08 & 3.70 & 1.09E-09 & 3.84 & 7.58E-11 \\
$Y({\rm N}_2)$          & 1.77E-07 & 3.95 & 1.15E-08 & 4.07 & 6.85E-10 \\
$\rho$                  & 5.00E-09 & 4.01 & 3.10E-10 & 4.09 & 1.82E-11 \\
$T$                     & 1.21E-02 & 4.02 & 7.44E-04 & 4.05 & 4.48E-05 \\
$\rho h$                & 6.77E+00 & 3.99 & 4.26E-01 & 4.07 & 2.54E-02 \\
\hline
	\end{tabular}
\end{center}	
\end{table}

\begin{figure}
\begin{center}
\includegraphics[width=5in]{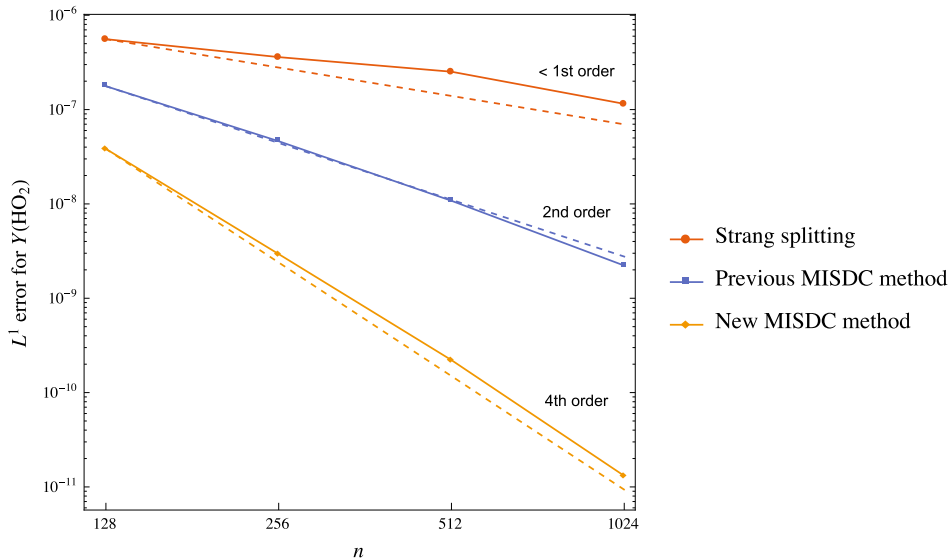}
\caption{\label{fig:algorithm_compare} Log-log plot of $L^1$ error 
as a function of resolution ($n$ is the number of cells)
for the intermediate species ${\rm HO}_2$ in a hydrogen flame using the
Strang splitting code \cite{DayBell:2000}, second-order MISDC algorithm \cite{LMC_SDC},
and the fourth-order MISDC algorithm presented in this paper.  Solid lines are
the numerical results, and dotted lines give first, second, and
fourth-order reference slopes.}
\end{center}
\end{figure}

We now compare the accuracy of our new code to our previous second-order
algorithm \cite{LMC_SDC} and the previous Strang splitting algorithm
\cite{DayBell:2000}.  In Figure \ref{fig:algorithm_compare}, we
plot the $L^1$ error for the intermediate species ${\rm HO}_2$
using each algorithm.  We use the problem setup described in Section
5.2 of \cite{LMC_SDC}, which is the same as described above except that
the time step is 20\% larger than described above at each resolution.
Note that not only is our new algorithm fourth-order, but at coarse resolution
($n=128$ with a domain length of 1.2~cm) the error in our fourth-order
code is already a factor of 5 smaller than it was with our previous second-order code.
In these simulations, the Strang split algorithm is not in the asymptotic
convergence regime, and exhibits only first-order behaviour, as noted
in \cite{LMC_SDC}.

We note that performing fewer than eight MISDC iterations in the fourth-order algorithm
causes an observed order-reduction in
our simulations unless we decrease the time step size as well.  For the cases analysed
when using four MISDC iterations, fourth-order is realised only when the 
time step is reduced by an order of magnitude or more.
This is not particularly surprising given the disparity in time
scales between the physical processes.  In these cases, we found it more efficient computationally
to increase the time step size, even though more MISDC iterations are required per step.

\begin{figure}
\begin{center}
\includegraphics[width=5in]{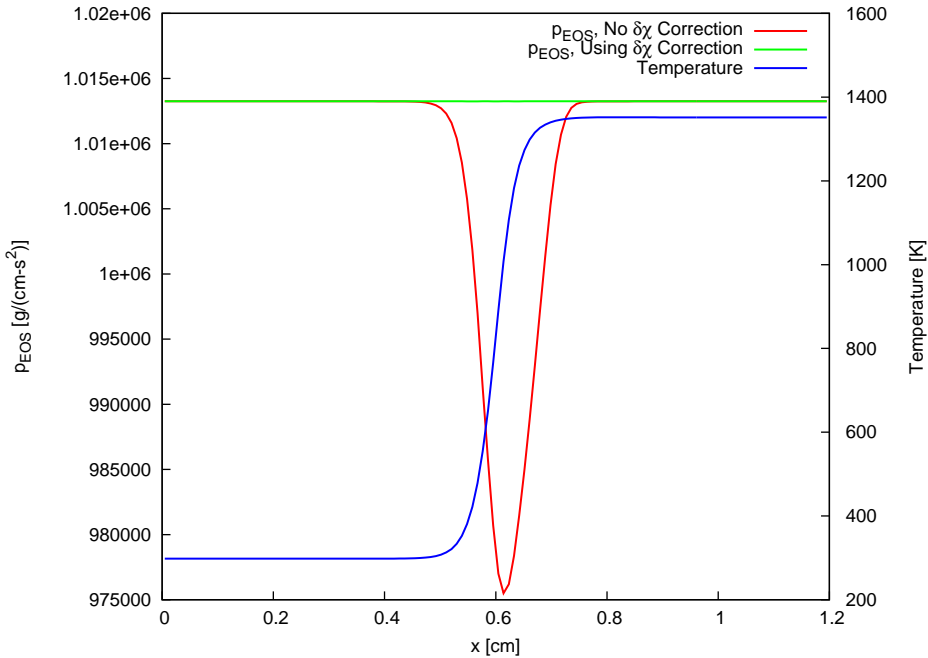}
\caption{\label{fig:deltachi_compare1}  A plot of the thermodynamic pressure, $p_{\rm EOS}$,
as a function of space for the $n=128$ simulation at the final time with and without the
$\delta\chi$ volume discrepancy correction.
The blue temperature plot is included for reference, indicating the location and shape of the flame.}
\end{center}
\end{figure}

\begin{figure}
\begin{center}
\includegraphics[width=5in]{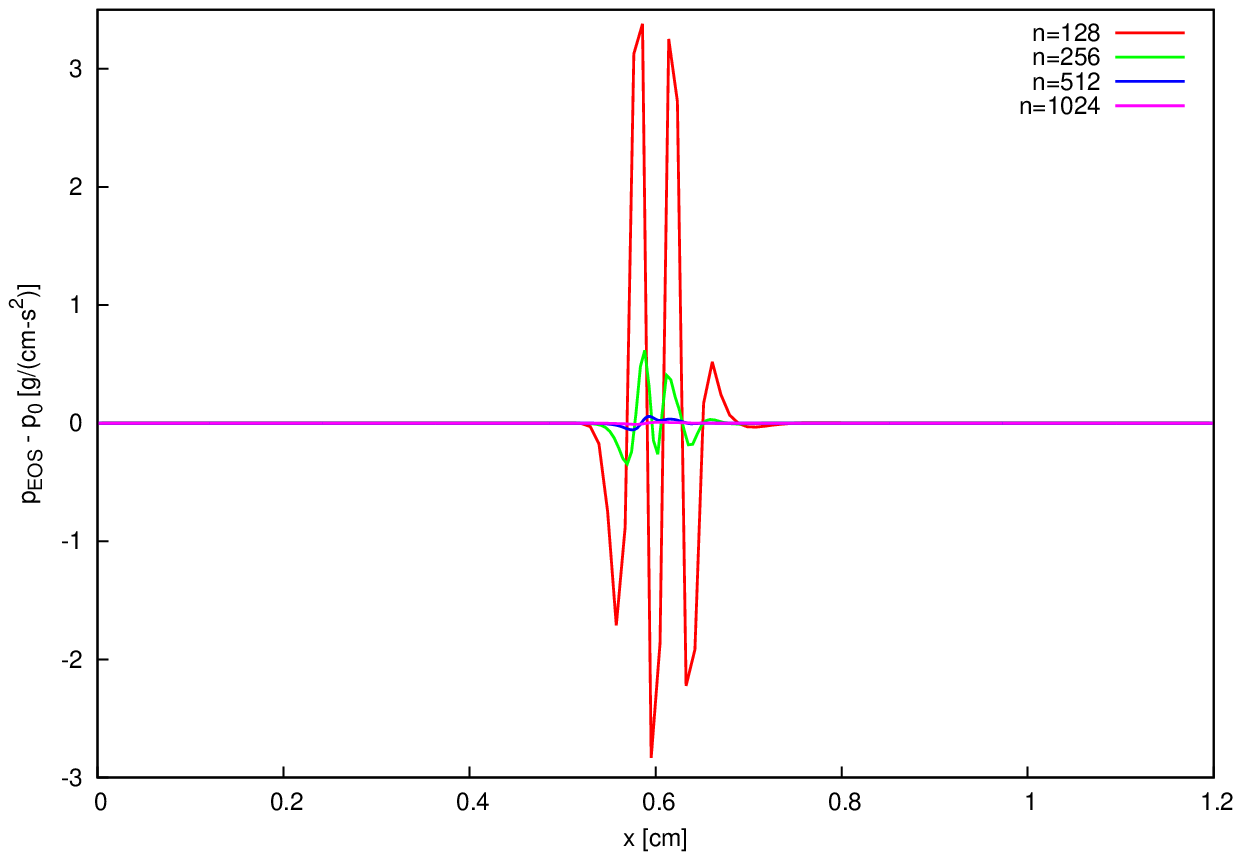}
\caption{\label{fig:deltachi_compare2} A plot of the thermodynamic drift, $p_{\rm EOS}-p_0$,
as a function of space for the $n=128, 256, 512,$ and $1024$ simulations with the
$\delta\chi$ correction term.}
\end{center}
\end{figure}

Next, we examine the effectiveness of the volume discrepancy algorithm in reducing the
thermodynamic drift.  We performed the exact same set of simulations described above,
but disabled the $\delta\chi$ volume discrepancy correction term.  In Figure 
\ref{fig:deltachi_compare1}, we plot the thermodynamic pressure, $p_{\rm EOS}$,
as a function of space for the $n=128$ simulation with and without the $\delta\chi$ correction
terms.  The correction causes the pressure to stay on the equation of state to within 
$\sim 3$~g/(cm-s$^2$),
whereas without the correction, the pressure drifts by $\sim 40,000$~g/(cm-s$^2$).
The blue temperature
plot is included for reference, indicating the location and shape of the flame.  Next,
in Figure \ref{fig:deltachi_compare2} we show a plot of the thermodynamic drift, $p_{\rm EOS}-p_0$,
as a function of space for the $n=128, 256, 512,$ and $1024$ simulations with the
$\delta\chi$ correction term.  The $\delta\chi$ correction term clearly drives
the pressure drift to zero as spatial resolution increases,
as the magnitude of the drift decreases by roughly a factor of 8 as
we increase resolution by a factor of 2.  For the $n=1024$ simulation, the maximum
value of $|p_{\rm EOS}-p_0|$ is less than 0.01~g/(cm-s$^2$).

\subsubsection{Methane Flame}
We next study the performance of the MISDC algorithm using a one-dimensional
premixed methane flame.  The simulations are based on the
GRIMech-3.0 \cite{GRI30:1995} model and associated database, as given in the
CHEMKIN-III library \cite{CHEMKIN:1996} format.  The GRIMech-3.0 model consists
of 53 species with a 325-step chemical reaction network for premixed methane
combustion.  This example is particularly
challenging because of the extremely broad range of
chemical time scales, which range from $10^{-4}$ to $10^{-10}$ seconds.

Similar to the hydrogen flame, the initial conditions
are obtained by interpolating from a frame-shifted, refined steady,
one-dimensional solution computed using the PREMIX code.  For this case, the
inlet stream at $T$=298~K and $p$=1~atm has composition, $Y({\rm O}_2 : {\rm CH}_4 : {\rm
N}_2) = (0.2238 : 0.0392 : 0.7370)$ so that the unstrained laminar burning
speed is $s_L = 18.9$~cm~s$^{-1}$.  The PREMIX solution is
interpolated onto a 1.2~cm domain with 8192 uniform cells, and
evolved with the second-order code from \cite{LMC_SDC} for $80~\mu$s with 
$\Delta t = 0.4~\mu$s.
The resulting solution is averaged down to coarse uniform meshes of 
$n =$ 128, 256, 512, and 1024 cells to provide the initial 
conditions for testing the fourth-order algorithm.  We evolve the system for 
$800~\mu$s to allow the solution to propagate $\sim$111~$\mu$m across the mesh
with $\sigma$$\approx$ 0.21.
The error and convergence results for the primary reactants, products, and key
intermediate species, as well as the remaining thermodynamic variables
are presented in Table \ref{table:GRI}.  We see 
fourth-order convergence in all variables, noting that for $Y({\rm CH}_3)$, 
the profile is extremely thin so
that higher resolution is required to reach the asymptotic regime.
\begin{table}[tb]
\caption{Error and convergence rates for a premixed methane flame using 
the fourth-order MISDC method with an advective CFL of $\sigma\approx 0.21$.
\label{table:GRI}}
\vspace{10pt}
\begin{center}
	\begin{tabular}{ l c c c c c }
	\hline
	Variable & $L^1_{128}$ & $r^{128/256}$ & $L^1_{256}$ & $r^{256/512}$ & $L^1_{512}$\\
	\hline\\[-1.5ex]
$Y({\rm CH}_4)$          & 1.11E-06 & 4.00 & 6.97E-08 & 3.98 & 4.42E-09 \\
$Y({\rm O}_2)$           & 3.77E-06 & 3.96 & 2.42E-07 & 4.07 & 1.44E-08 \\
$Y({\rm H}_2{\rm O})$    & 2.30E-06 & 4.02 & 1.42E-07 & 4.05 & 8.53E-09 \\
$Y({\rm CO}_2)$          & 1.87E-06 & 4.02 & 1.15E-07 & 4.07 & 6.87E-09 \\
$Y({\rm CH}_3)$          & 3.11E-08 & 2.48 & 5.59E-09 & 3.75 & 4.16E-10 \\
$Y({\rm CH}_2{\rm (S)})$ & 8.01E-11 & 4.14 & 4.54E-12 & 3.85 & 3.15E-13 \\
$Y({\rm O})$             & 1.05E-07 & 4.08 & 6.20E-09 & 3.90 & 4.16E-10 \\
$Y({\rm H})$             & 3.48E-09 & 3.83 & 2.45E-10 & 3.81 & 1.75E-11 \\
$Y({\rm N}_2)$           & 3.58E-07 & 3.74 & 2.68E-08 & 4.00 & 1.67E-09 \\
$\rho$                   & 1.25E-08 & 4.03 & 7.64E-10 & 4.05 & 4.61E-11 \\
$T$                      & 3.52E-02 & 4.01 & 2.18E-03 & 4.06 & 1.31E-04 \\
$\rho h$                 & 4.09E+01 & 3.97 & 2.60E+00 & 4.00 & 1.62E-01 \\
\hline
	\end{tabular}
\end{center}	
\end{table}

\subsubsection{Dimethyl Ether Flame}
Finally, we present a one-dimensional simulation of a premixed flame using a 39-species, 175-reaction
dimethyl ether (DME) chemistry mechanism \cite{Bansal:2015}.
The DME mechanism used in this test is extremely stiff.  It is quite challenging to 
capture the nonlinear coupling between diffusion and reaction chemistry, while evolving the
system on the much slower advection scales.
The inlet stream at $T$=298~K has composition, $Y({\rm CH}_3{\rm OCH}_3 : {\rm O}_2 : {\rm
N}_2) = (0.0726 : 0.2160 : 0.7114)$, and $p$=1~atm; the unstrained laminar burning
speed is $s_L = 24.9$~cm~s$^{-1}$.  The initial PREMIX-computed profiles
are interpolated onto a 0.6~cm domain with 8192 uniform cells and evolved for $16~\mu$s
using the second-order algorithm from \cite{LMC_SDC}
at $\Delta t = 0.08~\mu$s.
The resulting solution is averaged onto uniform grids 
of $n =$ 128, 256, 512, and 1024 cells to provide initial data for testing 
our fourth-order algorithm.  We evolve the system for 320~$\mu$s to
allow the solution to propagate 64~$\mu$m across the domain at $\sigma\approx$ 0.25.
The error and convergence  results for the primary reactants, products, and key intermediate species
\cite{Bansal:2015}, as well as the remaining thermodynamic variables
are presented in Table \ref{table:DME}.  We see fourth-order convergence in all variables.
\begin{table}[tb]
\label{table:DME}
\caption{Error and convergence rates for a dimethyl ether flame using 
the fourth-order MISDC method with an advective CFL of $\sigma\approx 0.25$.
	\label{table:DME}}
\vspace{10pt}
\begin{center}
	\begin{tabular}{ l c c c c c }
	\hline
	Variable & $L^1_{128}$ & $r^{128/256}$ & $L^1_{256}$ & $r^{256/512}$ & $L^1_{512}$\\
	\hline\\[-1.5ex]
$Y({\rm CH}_3{\rm OCH}_3)$           & 2.29E-06 & 3.83 & 1.62E-07 & 3.93 & 1.06E-08 \\
$Y({\rm O}_2)$                       & 2.99E-06 & 3.63 & 2.42E-07 & 4.02 & 1.49E-08 \\
$Y({\rm CO}_2)$                      & 2.51E-06 & 3.83 & 1.76E-07 & 4.04 & 1.07E-08 \\
$Y({\rm H}_2{\rm O})$                & 1.62E-06 & 3.51 & 1.42E-07 & 4.01 & 8.85E-09 \\
$Y({\rm CH}_3{\rm OCH}_2{\rm O}_2)$  & 1.55E-10 & 4.51 & 6.76E-12 & 3.88 & 4.61E-13 \\
$Y({\rm OH})$                        & 3.24E-07 & 3.80 & 2.32E-08 & 4.02 & 1.43E-09 \\
$Y({\rm HO}_2)$                      & 1.46E-07 & 3.80 & 1.05E-08 & 3.95 & 6.77E-10 \\
$Y({\rm O})$                         & 1.70E-07 & 3.55 & 1.46E-08 & 3.92 & 9.66E-10 \\
$Y({\rm H})$                         & 8.35E-09 & 3.68 & 6.52E-10 & 3.96 & 4.20E-11 \\
$Y({\rm N}_2)$                       & 1.09E-06 & 3.76 & 8.01E-08 & 3.93 & 5.25E-09 \\
$\rho$                               & 9.44E-09 & 3.67 & 7.42E-10 & 4.02 & 4.58E-11 \\
$T$                                  & 2.54E-02 & 3.59 & 2.11E-03 & 4.01 & 1.31E-04 \\
$\rho h$                             & 5.89E+01 & 3.83 & 4.15E+00 & 4.02 & 2.56E-01 \\
\hline
	\end{tabular}
\end{center}	
\end{table}

\makeatletter{}\section{Conclusions and Future Work}\label{sec:conclusions}
We have developed a fourth-order finite-volume algorithm for low Mach number
reacting flow with detailed kinetics and transport.  The approach iteratively couples
advection, diffusion, and reaction processes using efficient numerical methods
for each step.  The method exhibits much greater accuracy, even at coarse resolution,
than our previous second-order deferred correction strategy \cite{LMC_SDC} 
and Strang splitting algorithms \cite{DayBell:2000}.
We have incorporated a
volume discrepancy scheme that allows us to simultaneously conserve mass and enthalpy
while satisfying the equation of state to a high degree of accuracy.  The
volume discrepancy scheme is iterative, and naturally fits within our MISDC framework
with negligible computational cost.  We have
discussed an instability with our previous development path that did not allow
the method to extend to higher-order and demonstrated that our approach is
stable and convergent for a much broader range of parameters.

As discussed in the Results section, a key parameter in our new scheme
is the number of SDC iterations taken on each time step.  Formally,
the algorithm requires four iterations to couple together all
processes and achieve fourth-order convergence behavior.  However,
because our difficult demonstration problems feature cases with
relatively stiff diffusion and reaction processes, the inter-process
coupling at the advection time scale was not sufficiently accurate to
achieve the design rates.  There were at least two remedies to improve
the coupling: reduce the time step or increase the number of
iterations per step -- both of which make the algorithm more costly in
different ways.  By trial and error, we found that 8 iterations per
step was sufficient to achieve fourth-order accuracy for all cases
presented.  It is likely that an adaptive procedure can be developed
to optimize this choice for the general case.

The long-term goal of this effort is to extend the method described here 
to multidimensional, adaptive mesh refinement (AMR) simulations.  One 
issue concerning this goal is extending the projection
method formulation to multiple dimensions, where the velocity field is no longer
uniquely specified by the boundary conditions and the thermodynamic state.
Previous high-order SDC algorithms for incompressible 
flows (e.g. \cite{Kadioglu:2008,AABM:2011})
have been based on a gauge variable formulation which does not immediately extend to more general low Mach number models. 
One possible path forward is to utilise finite volume
Stokes solvers to allow us to solve the coupled viscous/projection step to
arbitrary spatial accuracy, and incorporate this into a method of lines approach
to allow for higher-order temporal integration.  Some work has already been done in
efficient projection-preconditioned finite volume Stokes solvers \cite{Cai:2014},
which fit well in our SDC based algorithms.
We would also like to implement an SDC based AMR algorithm that subcycles in time 
as in multilevel SDC methods \cite{Speck:MLSDC}.  
Finally, including stratification in our low Mach number reacting flow
models will allow them to be 
used in atmospheric \cite{Duarte:2015} and astrophysical \cite{MAESTRO} simulations.

\section*{Acknowledgements}
The work at LBNL was supported by the Applied Mathematics Program
of the DOE Office of Advanced Scientific Computing Research
under the U.S. Department of Energy under contract DE-AC02-05CH11231.

\makeatletter{}\appendix
\section{Convergence analysis of the previous method}\label{app:stability}
Here we examine the convergence properties of MISDC iterations of the method 
described in \cite{LMC_SDC}. As in Section \ref{subsec:misdc-stability}, we consider 
the linear ODE
\begin{align}
   \label{eq:linear-ode-stability}
   \phi_t = a\phi + d\phi + r\phi \equiv F(\phi),
\end{align}
For simplicity we will consider only two temporal nodes, $t^{n,0}=t^n$ 
and $t^{n,1} = t^{n+1} = t^n + \Delta t$. The corresponding Gauss-Lobatto quadrature 
rule in this case is the trapezoidal rule.
We compute a provisional solution using the diffusion correction equation 
using the implicit formula
\begin{eqnarray}
\label{eq:diff-corr-eqn-stability}
\phi^{1,(k+1)}_{\rm AD} = \phi^n + \Delta t^m
                        \left[d\phi_{\rm AD}^{1, (k+1)}  - d\phi^{1,(k)}\right] + I\left[(a+d)\phi^{(k)}\right] + I_R^{(k)}.
\end{eqnarray}
The term $I_R^{(k)}$ is equal to the integral of the reaction term.
The reaction correction equation is differentiated to obtain the ODE
\begin{equation}
   \phi^{(k+1)}_t(t) = r\phi^{(k+1)}(t) + d\phi_{\rm AD}^{1, (k+1)}  - d\phi^{1,(k)} 
                     + A(\phi^{(k)}, t) + D(\phi^{(k)}, t).
   \label{eq:reaction-corr-eqn-stability}
\end{equation}
Here $A$ and $D$ are polynomials representing the advection and diffusion contributions. We require 
that the integrals of $A$ and $D$ are equal to the numerical quadrature. In 
\cite{LMC_SDC}, $A$ and $D$ are chosen to be the averages at times $t^n$ and 
$t^n + \Delta t$. In order for this method to generalise to higher order, we 
would need to represent $A$ and $D$ by higher degree polynomials. For instance, 
we can take $A$ and $D$ to be the linear interpolants given by
\begin{align}
   A(\phi^{(k)}, t) &= (1-t/\Delta t)a\phi^{0,(k)} + (t/\Delta t)a\phi^{1,(k)}, \\
   D(\phi^{(k)}, t) &= (1-t/\Delta t)d\phi^{0,(k)} + (t/\Delta t)d\phi^{1,(k)},
\end{align}
so that their integrals are equal to the quadrature computed using the trapezoid 
rule. The reaction integral can then be computed by integrating both sides of 
(\ref{eq:reaction-corr-eqn-stability}), and rearranging:
\begin{eqnarray}
   I_R^{(k+1)} &\equiv& \int_{t^n}^{t^n + \Delta t} r\phi^{(k+1)}(\tau)\ d\tau\nonumber\\
&=& \phi^{1,(k+1)} - \phi^{0,(k+1)} 
             + \Delta t \left[\frac{a\phi^{0,(k)} + a\phi^{1,(k)}}{2} 
             + \frac{d\phi^{0,(k)} - d\phi^{1,(k)}}{2}
             + d\phi_{\rm AD}^{1, (k+1)} \right].\nonumber\\
\end{eqnarray}
We notice that equation (\ref{eq:reaction-corr-eqn-stability}) is an ODE of the form
\begin{align*}
   y_t(t) &= ry + c_1t + c_0,\\
   y(0) &= 0,
\end{align*}
to which the exact solution is given by
\begin{equation}
   \label{eq:ode-exact-soln}
   y(t) = -\frac{1}{r^2}\left(
      c_1 r t  - c_1 e^{rt} + c_0 r - c_0 e^{rt} + c1
   \right).
\end{equation}

Expanding expressions (\ref{eq:diff-corr-eqn-stability}) and 
(\ref{eq:reaction-corr-eqn-stability}), using the solution given by 
(\ref{eq:ode-exact-soln}), we see that the difference between successive 
iterates is given by
\begin{equation*}
   \phi^{1,(k+1)} - \phi^{1,(k)} = \alpha \left[ \phi^{1,(k)} - \phi^{1,(k-1)} \right]
                                 + \beta \left[ I_R^{(k)} - I_R^{(k-1)} \right],
\end{equation*}
where, in the case of the linear ODE (\ref{eq:linear-ode-stability}), $\alpha$ and $\beta$ 
are given by
\begin{align*}
   \alpha &= -\tfrac{a \left(2 d \Delta t+\Delta t r (d \Delta t-2)+e^{\Delta t r} (d \Delta t (\Delta t r-2)+2)-2\right)+d \left(e^{\Delta t r} \left(d \Delta t^2 r-2 \Delta t (d+r)+2\right)+d \Delta t (\Delta t r+2)-2\right)}{2 \Delta t r^2 (d \Delta t-1)}\\
   \beta &= \tfrac{d(e^{r\Delta t} - 1)}{r(d\Delta t - 1)}.
\end{align*}
We note that a sufficient condition for the iterative scheme to converge is 
$|\alpha|, |\beta| < 1$. For fixed $a$, we observe that the set of parameters $(d, r)$ 
that satisfy this condition (Figure \ref{fig:stability}) is exceedingly small.

Furthermore, making the ansatz that $a=\tilde a/\Delta x$, $d = \tilde d/\Delta x^2$, 
and $\Delta t = \lambda\Delta x$, we can compute the limits of $\alpha$ and 
$\beta$ as $\Delta x$ tends to zero. We see
\begin{align*}
   \lim_{\Delta x \to 0} \alpha = \frac{6-\tilde d\tilde r\lambda^2}{12}, \qquad \lim_{\Delta x \to 0} \beta = 1.
\end{align*}
We therefore conclude that the sufficient condition for the iterative scheme 
to converge is met only for a limited choice of coefficients $\tilde a$, $\tilde d$ and $\lambda$.
This is in contrast to the present method, for which convergence is guaranteed 
as $\Delta x$ tends to zero, for any choice of coefficients, as demonstrated 
in Section \ref{subsec:misdc-stability}. Indeed, numerical experiments indicate that this 
method suffers from extremely restrictive time-step conditions in order to obtain convergence.

\bibliographystyle{tCTM}
\bibliography{tCTMguide}

\newcommand{\noopsort}[1]{} \newcommand{\printfirst}[2]{#1}
  \newcommand{\singleletter}[1]{#1} \newcommand{\switchargs}[2]{#2#1}
\begin{thebibliography}{28}
\providecommand{\natexlab}[1]{#1}

\bibitem[1]{AABM:2011}
A.S. Almgren, A.J. Aspden, J.B. Bell, and M. Minion, {\itshape On the use of
  higher-order projection methods for incompressible turbulent flow}, SIAM J.
  Sci. Comput. 53 (2013).

\bibitem[2]{LMC_SDC}
A. Nonaka, J.B. Bell, M.S. Day, C. Gilet, A.S. Almgren, and M.L. Minion,
  {\itshape A deferred correction coupling strategy for low Mach number flow
  with complex chemistry}, Combust. Theory Modelling  (2012) available at
  http://www.tandfonline.com/doi/full/10.1080/13647830.2012.701019.

\bibitem[3]{Dutt:2000}
A. Dutt, L. Greengard, and V. Rokhlin, {\itshape Spectral Deferred Correction
  Methods for Ordinary Differential Equations}, BIT 40 (2000), pp. 241--266.

\bibitem[4]{Minion:2003}
M.L. Minion, {\itshape Semi-Implicit Spectral Deferred Correction Methods for
  Ordinary Differential Equations}, Comm. Math. Sci. 1 (2003), pp. 471--500.

\bibitem[5]{BLM:2003}
A. Bourlioux, A.T. Layton, and M.L. Minion, {\itshape High-Order Multi-Implicit
  Spectral Deferred Correction Methods for Problems of Reactive Flow}, Journal
  of Computational Physics 189 (2003), pp. 651--675.

\bibitem[6]{Layton:2004}
A.T. Layton and M.L. Minion, {\itshape Conservative Multi-Implicit Spectral
  Deferred Correction Methods for Reacting Gas Dynamics}, Journal of
  Computational Physics 194 (2004), pp. 697--715.

\bibitem[7]{zadunaisky:1964}
P. Zadunaisky, {\itshape A method for the estimation of errors propagated in
  the numerical solution of a system of ordinary differential equations}, in
  {\itshape The Theory of Orbits in the Solar System and in Stellar Systems.
  Proceedings of International Astronomical Union, Symposium 25}G.~Contopoulos
  ed.,  , Thessaloniki, 1964, pp. 281--287.

\bibitem[8]{pember-flame}
R.B. Pember, L.H. Howell, J.B. Bell, P. Colella, W.Y. Crutchfield, W.A.
  Fiveland, and J.P. Jessee, {\itshape An Adaptive Projection Method for
  Unsteady Low-{M}ach Number Combustion}, Comb. Sci. Tech. 140 (1998), pp.
  123--168.

\bibitem[9]{DayBell:2000}
M.S. Day and J.B. Bell, {\itshape Numerical simulation of laminar reacting
  flows with complex chemistry}, Combust. Theory Modelling 4 (2000), pp.
  535--556.

\bibitem[10]{rehmBaum:1978}
R.G. Rehm and H.R. Baum, {\itshape The equations of motion for thermally driven
  buoyant flows}, Journal of Research of the National Bureau of Standards 83
  (1978), pp. 297--308.

\bibitem[11]{majdaSethian:1985}
A. Majda and J.A. Sethian, {\itshape Derivation and Numerical Solution of the
  Equations of Low {M}ach Number Combustion}, Comb. Sci. Tech. 42 (1985), pp.
  185--205.

\bibitem[12]{Kee:1983}
R.J. Kee, J. Warnatz, and J. Miller, {\itshape Fortran computer-code package
  for the evaluation of gas-phase viscosities, conductivities, and diffusion
  coefficients.}, NTIS, SPRINGFIELD, VA(USA), 1983, 37  (1983).

\bibitem[13]{Warnatz:1982}
J. Warnatz, {\itshape Influence of transport models and boundary conditions on
  flame structure}, , in {\itshape Numerical methods in laminar flame
  propagation}   Springer, 1982, pp. 87--111.

\bibitem[14]{Layton:2005}
A.T. Layton and M.L. Minion, {\itshape Implications of the Choice of Quadrature
  Nodes for {P}icard Integral Deferred Corrections Methods for Ordinary
  Differential Equations}, BIT 45 (2005), pp. 341--373.

\bibitem[15]{Najm:1998}
H.N. Najm, P.S. Wyckoff, and O.M. Knio, {\itshape A semi-implicit numerical
  scheme for reacting flow. {I}. {S}tiff chemistry}, J. Comp. Phys. 143 (1998),
  pp. 381--402.

\bibitem[16]{Knio:1999}
O.M. Knio, H.N. Najm, and P.S. Wyckoff, {\itshape A semi-implicit numerical
  scheme for reacting flow. {I}. {S}tiff, operator-split formulation}, J. Comp.
  Phys. 154 (1999), pp. 428--467.

\bibitem[17]{Perini:2012}
F. Perini, E. Galligani, and R.D. Reitz, {\itshape An Analytical Jacobian
  Approach to Sparse Reaction Kinetics for Computationally Efficient Combustion
  Modeling with Large Reaction Mechanisms}, Energy \& Fuels 26 (2012), pp.
  4804--4822.

\bibitem[18]{Kadioglu:2008}
S.Y. Kadioglu, R. Klein, and M.L. Minion, {\itshape A Fourth-Order Auxiliary
  Variable Projection Method for Zero-{M}ach number gas dynamics}, Journal of
  Computational Physics 227 (2008), pp. 2012--2043.

\bibitem[19]{McCorquodale:2011}
P. McCorquodale and P. Colella, {\itshape A high-order finite-volume method for
  conservation laws on locally refined grids}, Communications in Applied
  Mathematics and Computational Science 6 (2011), pp. 1--25.

\bibitem[20]{Zhang:2012}
Q. Zhang, H. Johansen, and P. Colella, {\itshape A fourth-order accurate
  finite-volume method with structured adaptive mesh refinement for solving the
  advection-diffusion equation}, SIAM Journal on Scientific Computing 34
  (2012), pp. B179--B201.

\bibitem[21]{GRI30:1995}
M. Frenklach, H. Wang, M. Goldenberg, G.P. Smith, D.M. Golden, C.T. Bowman,
  R.K. Hanson, W.C. Gardiner, and V. Lissianski, {\itshape {GRI-Mech}---An
  Optimized Detailed Chemical Reaction Mechanism for Methane Combustion},
  GRI-95/0058, Gas Research Institute, 1995
  \url{http://www.me.berkeley.edu/gri_mech/}.

\bibitem[22]{CHEMKIN:1996}
R.J. Kee, R.M. Ruply, E. Meeks, and J.A. Miller, {\itshape CHEMKIN-{III}: A
  {FORTRAN} Chemical Kinetics Package for the Analysis of Gas-phase Chemical
  and Plasma Kinetics}, Technical Report SAND96-8216, Sandia National
  Laboratories, Livermore, 1996.

\bibitem[23]{PREMIX:1985}
R.J. Kee, J.F. Grcar, M.D. Smooke, and J.A. Miller, {\itshape {PREMIX}: A
  Fortran Program for Modeling Steady, Laminar, One-Dimensional Premixed
  Flames}, Technical Report SAND85-8240, Sandia National Laboratories,
  Livermore, 1983.

\bibitem[24]{Bansal:2015}
G. Bansal, A. Mascarenhas, and J.H. Chen, {\itshape Direct numerical
  simulations of autoignition in stratified dimethyl-ether (DME)/air turbulent
  mixtures}, Combustion and Flame 162 (2015), pp. 688 -- 702.

\bibitem[25]{Cai:2014}
M. Cai, A. Nonaka, J.B. Bell, B.E. Griffith, and A. Donev, {\itshape Efficient
  Variable-Coefficient Finite-Volume {S}tokes Solvers}, Commun. Comput. Phys.
  16 (2014), pp. 1263--1297.

\bibitem[26]{Speck:MLSDC}
R. Speck, D. Ruprecht, M. Emmett, M. Minion, M. Bolten, and R. Krause,
  {\itshape A multi-level spectral deferred correction method}, BIT Numerical
  Mathematics 55 (2015), pp. 843--867.

\bibitem[27]{Duarte:2015}
M. Duarte, A.S. Almgren, and J.B. Bell, {\itshape A low {M}ach number model for
  moist atmospheric flows}, J. Atmos. Sci. 72 (2014), pp. 1605--1647.

\bibitem[28]{MAESTRO}
A. Nonaka, A. Almgren, J. Bell, M. Lijewski, C. Malone, and M. Zingale,
  {\itshape MAESTRO: An adaptive low Mach number hydrodynamics algorithm for
  stellar flows}, The Astrophysical Journal Supplement Series 188 (2010), p.
  358.

\end{thebibliography}

\end{document}